\documentclass[11pt, a4paper,leqno]{article}
\usepackage[utf8]{inputenc}
\usepackage{amsmath}
\usepackage{amsfonts}
\usepackage{amssymb}
\usepackage{graphicx}
\usepackage{pifont}
\usepackage[left=2.5cm,right=2.5cm,top=3cm,bottom=3cm]{geometry}
\usepackage{enumerate}
\usepackage {shadow}
\usepackage{amsmath , amsfonts , amssymb , mathrsfs}
\usepackage[linkcolor=blue,colorlinks=true]{hyperref}
\usepackage[T1]{fontenc}
\usepackage{lmodern}
\usepackage[babel]{csquotes}
\usepackage{abstract,lipsum}
\usepackage{graphicx}
\usepackage{lmodern}
\usepackage{tocbibind}
\usepackage{microtype}
\usepackage{xcolor}
\usepackage{color}
\usepackage{fancybox}
\usepackage{varioref}
\usepackage{stmaryrd} 
\usepackage{textcomp}
\usepackage{pgf, tikz} 
\usepackage{setspace,booktabs}
\usepackage{enumitem}
\usepackage{parskip}
\usepackage{hyperref}
\usepackage{xpatch}  
\makeatletter   
\xpatchcmd{\@thm}{\thm@headpunct{.}}{\thm@headpunct{}}{}{}   
\DeclareSymbolFont{calletters}{OMS}{cmsy}{m}{n}
\DeclareSymbolFontAlphabet{\mathcal}{calletters}
\usepackage{fancyhdr}
\usepackage{fancyvrb}
\usepackage{datatool}
\usepackage{datatool-fp}
\usepackage{csvsimple}
\numberwithin{equation}{section}
\usepackage{float}
\usepackage{mathptmx}

\usepackage{xpatch}  
\makeatletter   
\xpatchcmd{\@thm}{\thm@headpunct{.}}{\thm@headpunct{}}{}{}   

\usepackage{fancyhdr}
\numberwithin{equation}{section}

\usepackage[T1]{fontenc}
\usepackage[utf8]{inputenc}

\newtheorem{thm}{Theorem}[section]

\newtheorem{lem}{Lemma}[section]

\newtheorem {ass}{Assumption}[section]

\newtheorem {prop}{Proposition}[section]
\newenvironment{pr}{{\bf{Proof :~}}} 

\newcommand{\ep}{\varepsilon}
\newcommand{\mt}{\mathsf{T}}
\newcommand{\mts}{\mathsf{T}_{\!S}}
\newcommand{\mti}{\mathsf{T}_{\!I}}
\newcommand{\mtr}{\mathsf{T}_{\!R}}
\newcommand{\mteps}{\mathsf{T}_{\!\ep,S}}
\newcommand{\mtepi}{\mathsf{T}_{\!\ep,I}}
\newcommand{\mtepr}{\mathsf{T}_{\!\ep,R}}
\newcommand{\mtepj}{\mathsf{T}_{\!\ep,J}}
\newcommand{\mtj}{\mathsf{T}_{\!J}}
\newcommand{\mtep}{\mathsf{T}_{\ep}}
\newcommand{\mcx}{\mathcal{X}_{\ep}}
\newcommand{\un}{\underset}

\title{A SIR model on a refining spatial grid\\ - Law of Large Numbers }
\author{ M. N'zi \thanks{Univ Félix Houphouët Boigny, modestenzi@yahoo.fr}  \and E. Pardoux \thanks{Aix Marseille Univ, CNRS, Centrale Marseille, I2M, Marseille, France, etienne.pardoux@univ-amu.fr }
	\and T. Yeo \thanks{Aix Marseille Univ, CNRS, Centrale Marseille, I2M, Marseille, France, tenan.yeo@univ-amu.fr}}

\begin{document}
	
	\maketitle
	\labelformat{rmq}{Remark~#1}
	\labelformat{prop}{Proposition~#1}
	\labelformat{lem}{Lemma~#1}
	\labelformat{ass}{Assumption~#1}
	\labelformat{thm}{Theorem~#1}
	\newcommand{\dis}{\displaystyle}
	\newcommand{\tep}{T_{\varepsilon}}
	\newcommand{\eps} {\varepsilon}
	\newcommand{\sep}{\mathcal{S}_{\varepsilon}}
	\newcommand{\iep}{\mathcal{I}_{\varepsilon}}
	\newcommand{\rep}{\mathcal{R}_{\varepsilon}}
	\newcommand{\T}  {\mathbb{T}}
	\newcommand{\E}  {\mathbb{E}}
	\newcommand{\D}  {\mathbb{D}}
	\renewcommand{\P}  {\mathbb{P}}
	\newcommand{\R}  {\mathbb{R}}
	\newcommand{\bV}  {\big\Vert}
	\newcommand{\BV}  {\Big\Vert}
	\newcommand{\nn}  {\noindent}
	\newcommand{\wtd}{\widetilde{\Delta}_{\varepsilon}}
	\newcommand{\fpr}  { \vspace{-0.7cm} \begin{flushright}
			$\blacksquare$
	\end{flushright}}
	\setlength{\parindent}{0pt} 
	
	\begin{abstract}
		\nn We study in this paper a compartmental SIR model for a population distributed in a bounded domain D of $\mathbb{R}^d$, d= 1, 2 or 3. We describe a spatial model for the spread of a disease on a grid of D. We prove two laws of large numbers. On the one hand, we prove that  the stochastic model converges to the  corresponding deterministic patch model as the size of the  population tends to infinity. On the other hand, by letting both the size of the population tend to infinity and the mesh of the grid go to zero, we obtain a law of large numbers in the supremum norm, where the limit is a diffusion SIR model in D.
	\end{abstract}
	
	Funding: T\'enan Yeo was supported by a thesis scholarship from the government of Ivory Coast, and
	a salary as instructor at University of Aix--Marseille, and the two other authors by their respective university.
	
	Conflict of Interest: The authors declare that they have no conflict of interest.

	
	\setcounter{section}{-1}
	
	\section{Introduction}
	\rhead{Introduction}
	
	There is by now a good number of books and a huge number of papers treating mathematical models of epidemics. Most of them treat deterministic models, while some of them discuss as well stochastic models. Let us quote among many others  Kermack $\& $ McKendrick (1927), Anderson $\& $ Britton (2000), Britton~$\& $~Pardoux (2019). These last two works show that the standard deterministic models are law of large numbers limits of individual--based stochastic models. They also study fluctuations around the law of large numbers limit, via the central limit theorem, and concerning the last reference, the large deviations. Those fluctuations allow to explain extinction of an endemic disease, which is a stable equilibrium of the deterministic model.
	
	The classical SIR model ignores the fact that a population spreads over a spatial region. However environmental heterogeneity, spatial connectivity and movement of individuals play important roles in the spread of infectious diseases. Spatially uniform models are not sufficient to give a realistic picture of the spread of the disease. There is by now quite an important literature on spatial epidemics model, both in discrete and in continuous space, see e.g. Allen, Bolker, Lou  $\&$ Nevai (2007) and Allen, Bolker, Lou  $\&$ Nevai~(2008), and the references therein.
	
	In the present paper, we consider both deterministic and stochastic models in discrete and continuous space. More precisely, we start with 
	an individual based stochastic model for a population with constant size $\mathbf{N}\ep^{-d}$, distributed on the nodes of a regular grid discretizing $[0,1]^d$, with $d=1,2$ or $3$ (we shall concentrate mainly on the case $d=2$, which seems to us most relevant).  Letting first $\mathbf{N}\to\infty$, while $\eps$, the mesh size, is kept fixed, we shall obtain as law of large numbers limit a system of ODEs on the grid, which is a patch epidemics model. Letting then $\eps\to0$, we will show that the system of ODEs converges to a system of PDEs on $[0,1]^d$, which is a deterministic epidemic model in continuous space. It is rather clear that one cannot hope to get the same result by letting first $\eps\to0$, and then $\mathbf{N}\to\infty$. Indeed, the first limit should be a continuous space model for quantities which take their values in the set $\{k/\mathbf{N},\, 0\le k\le\mathbf{N}\}$, with a partial differential operator for the displacement of the population, which would not make much sense. Consequently, if one wants to obtain a limit while letting jointly $\mathbf{N}\to\infty$ and $\eps\to0$, there must be a constraint which limits the speed of convergence of $\eps$ to $0$, in terms of the speed of convergence of $\mathbf{N}$ to $+\infty$. The weakest possible such constraint seems to be the one which has been first introduced by Blount (1992) for chemical reaction models, namely the restriction that $\mathbf{N}/\log(1/\eps)\to\infty$, see also Debussche $\&$ Nankep (2017). We shall extend that result to our situation where the limit is not a single PDE, but a system of PDEs.

	The model is constructed on a $d$--dimensional bounded domain $[0,1]^d$ $(d= 1, 2,3)$. We first suppose that the population is spatially distributed on the nodes of a grid $ D_{\varepsilon} := [0,1]^d \cap \varepsilon \mathbb{Z}^d=\left\{ x_i,\; 1\le i\le \varepsilon^{-d} \right\} $ of $[0,1]^d$, where $0<\varepsilon<1$ (two neighboring sites are at distance $ \varepsilon$ apart, see Figure 1). Nodes represent   communities in which the disease can  grow. The population is divided  in three compartments S, I and R. For a space-time coordinate (t, $x_i$), we denote by
	
	\vspace{-0.3cm}
	
	\begin{enumerate}
		\item[$\bullet$] $S^{\varepsilon}(t,x_i)$ the number of susceptibles at  site $x_i$ at time $t$,
		
		\vspace{-0.3cm}
		
		\item[$\bullet$] $I^{\varepsilon}(t,x_i)$ the number of infected at site $x_i$ at time $t$,
		
		\vspace{-0.3cm}
		
		\item[$\bullet$] $R^{\varepsilon}(t,x_i)$ the number of removed at site $x_i$ at time $t$.
	\end{enumerate}
	
	\vspace{-0.3cm}
	
	In this case the deterministic model is given by a system of ordinary differential equation (ODE) and the stochastic one by a jump Markov process.
	Note that  Arnold $\&$ Theodosopulu (1980), Kotelenez (1986), Blount (1992), and also some of the references therein, describe such spatial models for chemical reactions. The resulting process has one component and is compared with the corresponding deterministic model. 
	
	In the present paper, we focus  our attention on the law of large numbers. In future works, we intend to discuss the fluctuations around the law of large numbers.
	
	Let us briefly describe the content of this paper. In section 1, we introduce a deterministic model on the grid $ D_{\varepsilon}$ of the  bounded domain  $[0,1]^d$ and we recall the relation between this model and the limiting PDE model on $[0,1]^d$ as $ \varepsilon\to 0 $. Then we introduce the  stochastic model on the same grid for a population of total size $\mathbf{N}\ep^{-d}$. In section 2, we fix the parameter $ \varepsilon $ and let the  initial average number $ \mathbf{N}$ of individuals in each site tend to infinity: the limiting law of large numbers limit is the already introduced deterministic model.
	As $\varepsilon\to 0 $ our system of ODEs converges towards a system of PDEs. 
	Finally in  section 3, we prove a law of large numbers in the supremum norm when we let both the size of the population go to infinity and the mesh  of the grid go to zero, under the weak restriction that $\dfrac{\mathbf{N}}{\log(1/\varepsilon)}\longrightarrow \infty$. 
	\begin{center}
		\begin{tikzpicture}[x=1.3cm, y=1.3cm]
		\foreach \i in {0,...,6}{
			\draw (\i,0)--(\i,6);
			\foreach \j in {0,...,6}{
				\draw (0,\j)--(6,\j);
				\fill (\i,\j) circle(5pt);
			}
		}
		\end{tikzpicture} \\
		Figure 1- $[0,1]\times [0,1]$ grid
	\end{center}
	\section{The models}
	\rhead{The models}
	
	Suppose that individuals are living in the bounded domain 
	$D:=(0,1)^d\subset \mathbb{R}^d$.
	We consider an infectious disease which spreads in the population. 
	Consider at each point of a grid (see Figure 1 ) on the d-dimensional domain $D$  a deterministic and a stochastic SIR model,
	with migration between neighboring sites (two neighboring sites are at distance $ \varepsilon $ apart). We assume that the mesh size of the grid $ \varepsilon $ is such that  $\varepsilon^{-1}\in \mathbb{N} $, where $ \mathbb{N}$ is the set of positive integers. We assume that the studied epidemic concerns a population of fixed size. In this model, infections are local. We let $ \beta \,  :  \, \mathbb{R}^d\longrightarrow  \mathbb{R}_+$ and $ \alpha\, : \,  \mathbb{R}^d \longrightarrow  \mathbb{R}_+$ be continuous functions and we set   $ \bar{\beta}=\un{x\in D}{\sup}\beta(x)$ and $ \bar{\alpha}=\un{x\in D}{\sup}\alpha(x)$. For each site $x_i$ 
	
	\vspace{-0.5cm}
	
	\begin{enumerate}
		\item[$ \bullet$]  Susceptible individuals become infectious at rate $ \beta(x_i)\dfrac{S^{\varepsilon}(t,x_i)}{S^{\varepsilon}(t,x_i)+I^{\varepsilon}(t,x_i)+R^{\varepsilon}(t,x_i)}I^{\varepsilon}(t,x_i)$.\\ 
		Note that an individual chosen uniformly at random  site $x_i$ at time $t$ is susceptible with probability $\dfrac{S^{\varepsilon}(t,x_i)}{S^{\varepsilon}(t,x_i)+I^{\varepsilon}(t,x_i)+R^{\varepsilon}(t,x_i)}$;
		\item[$ \bullet$] each infectious recovers at rate $ \alpha(x_i) $, so the total recovery rate at time t is $\alpha(x_i)I^{\varepsilon}(t,x_i)$;
		
		\vspace{-0.3cm}
		
		\item[$ \bullet$] the migrations of susceptible, infected and removed individuals between location $x_i$ and its neighboring sites occur at rate $\dfrac{\mu_S}{\varepsilon^2}S^{\varepsilon}(t,x_i)$, $\dfrac{\mu_I}{\varepsilon^2}I^{\varepsilon}(t,x_i)$ and $\dfrac{\mu_R}{\varepsilon^2}R^{\varepsilon}(t,x_i)$ respectively. $ \mu_S$, $\mu_I$  and  $\mu_R$ are positive diffusion coefficients for the susceptible, infected and removed subpopulations, respectively.
	\end{enumerate}
	Here, we assume that the compartment $R$ contains individuals who are dead or who have recovered and have permanent immunity. We can assume boundary conditions of the Neumann or periodic type. In this paper, we focus  our attention on Neumann boundary conditions (representing a closed environment i.e. there is no flux of individuals through the boundary). The choice $D=(0,1)^d$  as the spatial domain is made for the sake of simplifying the analysis, but our results can be extended to  any bounded domain $ D \subset \mathbb{R}^d $, with a reasonably smooth boundary.

	Initially $\mathbf{N}\ep^{-d}$ individuals are distributed on  the grid. That is, there is an average of $\mathbf{N}$ individuals on each site. 
	We first introduce the deterministic model and then we construct the corresponding stochastic model. \\
	In the following  we use the generic notation $C$ for a positive constant, the value of which  may change from line to line. These constants can depend upon some  parameters of the model, as long as these are independent of  $ \varepsilon $ and $ \mathbf{N}$.
	
	\subsection{The deterministic model}
	\rhead{The deterministic model}
	
	The space is the grid $D_{\varepsilon}$ of D. In order to take into account Neumann boundary conditions, we add some fictitious sites which extend the grid outside the domain, as shown in  Figure 2 below. We denote by $ \partial_{\vec{n}.out} D_{\varepsilon}$ the set of those fictitious sites.
	We use the notation $ y_i \sim x_i $ to mean that the sites $y_i$ and $x_i$ are neighbors. Each interior point of $ D_{\varepsilon}$ has $2d$ neighbours. Each boundary point has at least one fictitious site among its neighbors.
	\begin{center}
		\hspace{1.2cm}	\begin{tikzpicture}[x=1.3cm, y=1.3cm]
		\draw[black!20,line width=1.1cm+5pt] (0.5,0.5)--(0.5,6.5)--(6.5,6.5)--(6.5,0.5)--cycle;
		\foreach \i in {1,...,6}{
			\draw (\i,1)--(\i,6);
			\foreach \j in {1,...,6}{
				\draw (1,\j)--(6,\j);
				\fill (\i,\j) circle(5pt);
			}
		}
		
		\foreach \i in {1,...,6}{
			\draw[dashed] (\i,0)--(\i,7);
			\foreach \j in {0,7}{
				\fill[red] (\i,\j) circle(5pt);
			}
		}
		\foreach \j in {1,...,6}{
			\draw[dashed] (0,\j)--(7,\j);
			\foreach \i in {0,7}{
				\fill[red] (\i,\j) circle(5pt);
				
			}
		}
		\draw (8.7,5) node{\begin{huge}\color{red}$\bullet$\end{huge}= fictitious sites} ;
		\draw (9.9,4) node{{\boldmath$ \partial_{\vec{n}.out} D_{\varepsilon}$}:= the set of  fictitious sites};
		\fill[blue] (4,4) node[below left] {$x_i$};
		
		\draw [blue, thick,>=latex,->] (3+.1,4+.05) to[out=15,in=165](4-.1,4+.05);
		\draw [blue, thick,>=latex,->] (4+.1,4+.05) to[out=15,in=165](5-.1,4+.05);

		\draw[blue, thick,>=latex,->] (4-.1,4-.05) to[out=-165,in=-15](3+.1,4-.05);
		\draw[blue, thick,>=latex,->] (5-.1,4-.05) to[out=-165,in=-15](4+.1,4-.05);
		
		\draw[blue, thick,>=latex,->] (4+.05,3+.1) to[out=75,in=-75](4+.05,4-.05);
		\draw[blue, thick,>=latex,->] (4+.05,4+.1) to[out=75,in=-75](4+.05,5-.1);
		
		\draw[blue, thick,>=latex,->] (4-.05,4-.1) to[out=-105,in=105](4-.05,3+.1);
		\draw[blue, thick,>=latex,->] (4-.05,5-.1) to[out=-105,in=105](4-.05,4+.1);
		\end{tikzpicture}\\
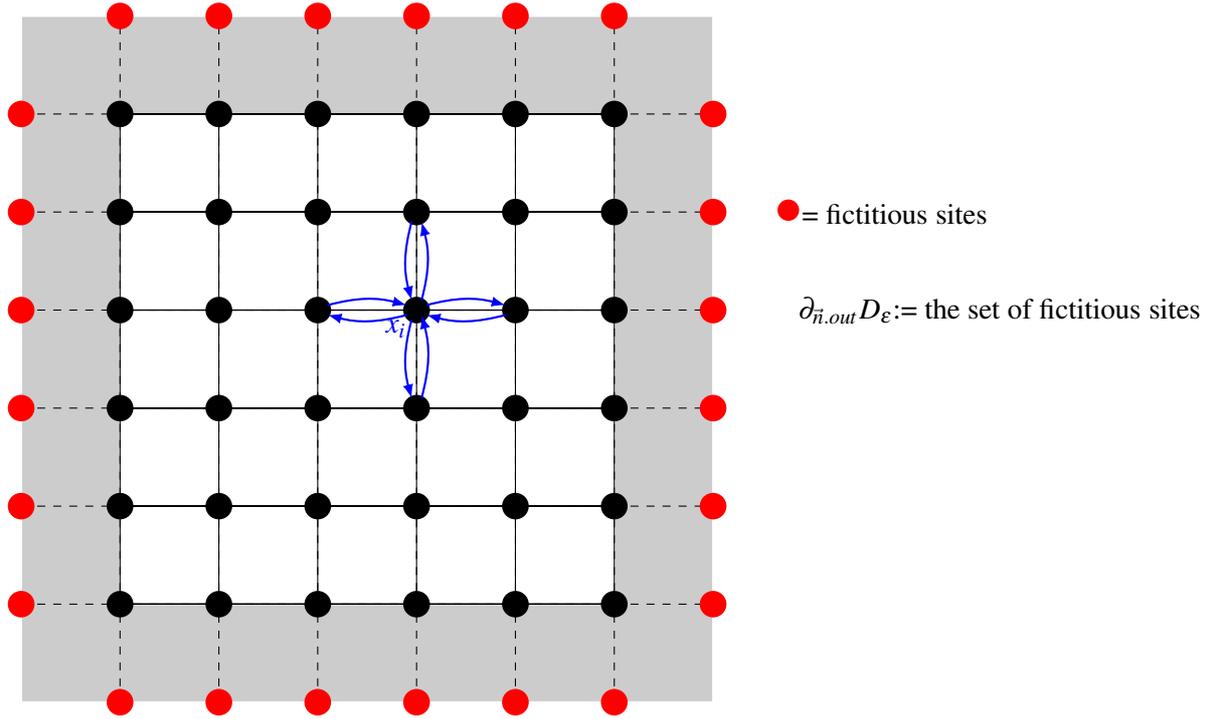
		$\hspace{-6cm} {\boldmath \text{Figure 2}-\text{Modeling the Neumann condition}}$
	\end{center}
	By thinking of an infinite size population allowing "proportions" in each compartment to be continuous, we have the following deterministic model for "proportions" (this point  of view will become quite clear in section 2 below):
	\begin{equation}\label{eqdet} 
	\hspace{-2cm}\left\{ 
	\begin{aligned} 
	\dfrac{d\,S_{\varepsilon}}{dt}(t,x_i) &=  - \dfrac{\beta(x_i)\, S_{\varepsilon}(t,x_i)I_{\varepsilon}(t,x_i)}{S_{\varepsilon}(t,x_i)+I_{\varepsilon}(t,x_i)+R_{\varepsilon}(t,x_i)}+ \mu_S\,\Delta_{\varepsilon} S_{\varepsilon}(t,x_i) \\
	\bigskip 
	\dfrac{d\,I_{\varepsilon}}{dt}(t,x_i) &=   \dfrac{\beta(x_i)\, S_{\varepsilon}(t,x_i)I_{\varepsilon}(t,x_i)}{S_{\varepsilon}(t,x_i)+I_{\varepsilon}(t,x_i)+R_{\varepsilon}(t,x_i)}-\alpha(x_i) \,I_{\varepsilon}(t,x_i) + \mu_I\,\Delta_{\varepsilon} I_{\varepsilon}(t,x_i) \\
	\dfrac{d\,R_{\varepsilon}}{dt}(t,x_i) & = \alpha(x_i) \,I_{\varepsilon}(t,x_i) +\mu_R\,\Delta_{\varepsilon} R_{\varepsilon}(t,x_i) , \;  \;
	(t,x_i)  \in (0,T)\times  D_{\varepsilon} \\
	& \left. 
	\hspace{-1.8cm} \begin{array}{rl}
	S_{\varepsilon}(t,x_i)= S_{\varepsilon}(t,y_i)\\
	I_{\varepsilon}(t,x_i)= I_{\varepsilon}(t,y_i)\\
	R_{\varepsilon}(t,x_i)= R_{\varepsilon}(t,y_i) 
	\end{array} \right\}  \text{for} \; x_i\in \partial D_{\varepsilon}, \;  x_i \sim y_i \; \text{and} \;  y_i \in \partial_{\vec{n}.out} D_{\varepsilon} \\
	& \hspace{-1.5cm}  S_{\varepsilon}(0,x_i), I_{\varepsilon}(0,x_i), R_{\varepsilon}(0,x_i)\ge 0, \; 0<S_{\varepsilon}(0,x_i)+ I_{\varepsilon}(0,x_i)+ R_{\varepsilon}(0,x_i) \le M,\\
	&\hspace{-1.5cm} \text{for some} \; M < \infty ,
	\end{aligned}
	\right. 
	\end{equation}
	where $ S_{\varepsilon}(t,x_i)$ (resp. $I_{\varepsilon}(t,x_i)$, resp. $R_{\varepsilon}(t,x_i)$) is  the proportion of the total population which is  both susceptible (resp. infectious, resp. removed) and located at site $x_i$ at time $t$. $\Delta_{\varepsilon}$ is the discrete Laplace operator defined as follows: 
	$ \dis\Delta_{\varepsilon}f(x_i) =  \varepsilon^{-2}\sum_{j=1}^{d}\big[f(x_i+\varepsilon e_j)-2f(x_i)+f(x_i-\varepsilon e_j) \big].$
	
	Note that (\ref{eqdet}) is the discrete space approximation  of the following system of PDE
	\begin{equation}\label{cdm}
	\hspace{-4.5cm}\left \{
	\begin{aligned}
	\dfrac{\partial\,\mathbf{s}}{\partial t}(t,x)= & \;-\dfrac{\beta(x)\, \mathbf{s}(t,x)\mathbf{i}(t,x) }{\mathbf{s}(t,x)+\mathbf{i}(t,x)+\mathbf{r}(t,x)}  + \mu_S \,\Delta \mathbf{s}(t,x) \\
	\dfrac{\partial\,\mathbf{i}}{\partial t}(t,x)=&\;\dfrac{\beta(x)\, \mathbf{s}(t,x)\mathbf{i}(t,x) }{\mathbf{s}(t,x)+\mathbf{i}(t,x)+\mathbf{r}(t,x)} - \alpha(x)\, \mathbf{i}(t,x)+\mu_I\, \Delta \mathbf{i}(t,x)\\
	\dfrac{\partial\,\mathbf{r}}{\partial t}(t,x)=&\;\alpha(x)\, \mathbf{i}(t,x)+\mu_S\, \Delta \mathbf{r}(t,x),  \quad (t,x) \in (0,T)\times  D \\
	&\hspace{-1.7cm} \dfrac{\partial \,\mathbf{s}}{\partial n_{\text{out}}}(t,x)=\;  \dfrac{\partial\, \mathbf{i}}{\partial n_{\text{out}}}(t,x)=\dfrac{\partial\, \mathbf{r}}{\partial n_{\text{out}}}(t,x) = 0, \; \;  \; \text{for}\; x\in \partial D \\
	&\hspace{-1.5cm}  \mathbf{s}(0,x), \mathbf{i}(0,x), \mathbf{r}(0,x)\ge 0 , \; 0< \mathbf{s}(0,x)+\mathbf{i}(0,x)+\mathbf{r}(0,x)\le M,
	\end{aligned}
	\right.
	\end{equation} 
	where $ \dfrac{\partial}{\partial n_{\text{out}}}$  denotes differentiation in the direction of the outward normal to $ \partial D$ and $ \Delta $ denotes the d-dimensional Laplace operator.\\ 
	System (\ref{cdm}) is a  reaction-diffusion epidemic model which has been studied by several authors. Webb~(1981) gave a similar reaction-diffusion model for a deterministic diffusive epidemic model, established the existence of solutions  and  analyzed their behavior as $ t\to \infty$. His method exploits tools of functional analysis and dynamical systems, specifically the theory of semigroups of linear and nonlinear operators in Banach spaces and Lyapunov stability techniques for dynamical systems in metric spaces. In the same way  Yamazaki $\&$  Wang (2016) gave  a reaction-convection-diffusion epidemic model for cholera dynamics and studied the global well-posedness and the asymptotic behavior of the solutions. See also Du $\&$ Peng (2016), Yamazaki~(2018a), Yamazaki (2018b). Let us mention that the SIR model (\ref{eqdet}) describes the spread of an infectious disease where recovered individuals gain immunity from re-infection. Of course in some cases recovered individuals have not permanent immunity. Hence individuals in the compartiment R  can experience reinfection. Moreover, susceptible individuals that become infected  can first pass through a latent stage (exposed). Such  models are used to study the transmission dynamics of the Ebola virus disease as treated in Agusto (2017a). Also, in Agusto et al. (2017b) the authors used  such model  to explore the Zika
	virus transmission dynamics in a human population.  Another  model which received  attention in the literature is the diffusion epidemic SIS model. In this model, when an infectious individual cures, he immediately becomes susceptible again. Such model has been considered in Allen et al. (2008). Although we restrict ourselves  to the SIR model, our results can easily  be adapted to SIRS, SIS, SEIR, SEIRS  models.
	
	Before describing the stochastic model, we introduce some notations and preliminaries, and then discuss the relation between the system of PDEs (\ref{cdm}) and its discretisation.
	
	\subsubsection{Some notations and preliminaries}
	
	In this subsection we introduce some notations and also give preliminary lemmas which will be  needed in our subsequent work. For all $ x_i \in D_{\varepsilon}$, let $V_i$ be the cube centered  at the site $ x_i$ with volume $ \varepsilon^d$.
	Let $ H^{\varepsilon} \subset L^2\big(D\big)$ denote the space of real valued step functions that are constant on each  cell $ V_i$. For $f\in H^{\varepsilon} $, let us  define
	\[ \nabla_{\varepsilon}^{j,+}f(x_i) =\frac{f(x_i +\varepsilon e_j)-f(x_i)}{\varepsilon},\]
	\[ \nabla_{\varepsilon}^{j,-}f(x_i) =\frac{f(x_i)- f(x_i-\varepsilon e_j)}{\varepsilon}. \] 
	It is not hard to see that
	$$ \langle \; \nabla_{\varepsilon}^{j,+}f, g \; \rangle = -\langle \; f  , \nabla_{\varepsilon}^{j,-}g \; \rangle , $$
	$$\Delta_{\varepsilon}f(x_i)=\sum_{j=1}^{d}\nabla_{\varepsilon}^{j,-}\nabla_{\varepsilon}^{j,+}f(x_i).$$

	We introduce the canonical projection $ \dis P_{\varepsilon} : L^2(D)\longrightarrow H^{\varepsilon} $ given by
	$$\varphi\longmapsto P_{\varepsilon}\varphi(x)= \varepsilon^{-d}\int_{V_i} \varphi(y)dy \; \; \; \; \; \text{if}\; x\in V_i . $$
	
	Throughout this paper, we assume that the initial condition satisfies
	\begin{ass}\label{a1}
		$ \mathcal{S}_{\varepsilon}(0,x)=P_{\varepsilon}\, \mathbf{s}(0,x)$, $ \mathcal{I}_{\varepsilon}(0,x)=P_{\varepsilon}\,\mathbf{i}(0,x)$, $ \mathcal{R}_{\varepsilon}(0,x)=P_{\varepsilon}\,\mathbf{r}(0,x)$  and 
		
		$ ~\hspace{2.7cm}\dis\int_D\big(\mathbf{s}(0,x)+\mathbf{i}(0,x)+\mathbf{r}(0,x)\big)dx=1$.
	\end{ass} 
	
	Here, we describe some of the spectral properties of the (discrete)-Laplacian  which will play an important  role in the sequel. More details can be found in Kotelenez(1986). 
	
	$\bullet $ For a multiindex  $ \dis m = (m_1 , \ldots , m _d ) $, where $m_j\in \mathbb{N}\cup \{0\}$, and $ x\in \mathbb{R} $, we define
	$$
	f_{m_j}(x)=\left\{
	\begin{array}{cc}
	& \hspace{-.5cm}\sqrt{2}\cos(m_j\pi  x), \; \text{for} \; \; m_j\ge 1 \\
	& 1 \; , \qquad \qquad  \; \;\, \mbox{for} \; \;  m_j= 0 \; .
	\end{array}
	\right.
	$$
	\newcommand{\f}{\mathbf{f}}
	For $ \varphi, \psi \in L^2\big(D\big)$, $\dis \langle \; \varphi , \phi \;  \rangle := \int_D \varphi(r)\phi(r)dr $  denotes the scalar product in $L^2\big(D\big).$\\
	For each $m\in \mathbb{Z}_+^d$, $x\!=\!(x^1,\cdots,x^d)\in D$, we define $\dis \f_{m}(x)\!=\!\prod_{j=1}^{d}f_{m_j}(x^j)$.   $\dis \big\{\, \f_{m},\, m\in \mathbb{Z}_+^d \,\big\} $ is a complete\\
	
	\vspace{-0.4cm}
	
	orthonormal system (CONS) of eigenvectors of $ \Delta $ in $L^2(D)$ with eigenvalues  $\dis -\lambda_{m}\!\!=\!\!-\pi^2\sum_{j=1}^{d}m_j^2. $
	Consequently, the semigroup $\mathsf{T}(t)\!:= \!\exp\big(\Delta \,t\big)$ acting on $ L^2\big(D\big)$  generated by $ \Delta $ can be represented by 
	
	\vspace{0.2cm}
	
	$\dis \hspace{3cm} \mathsf{T}(t)\varphi\!=\!\sum_m\exp(-\lambda_{m} t )\langle \; \varphi , \f_m \;  \rangle \f_m ,\; \; \varphi\in~
	L^2\big(D\big)$.
	
	\vspace{-0.2cm}
	
	$\bullet$ For $ i = ( i_1, \ldots , i_d)\in \left\{ 0, 1 , \ldots , \varepsilon^{-1}-1\right\}^d$,  let $\dis  	V_i = \prod_{j=1}^{d}\Big[\big(i_j-\dfrac{1}{2}\big)\varepsilon , \big(i_j+\dfrac{1}{2}\big)\varepsilon\Big) \subset [0,1]^d   $ 
	and for 
	
	\vspace{-0.2cm}
	
	$ m \in \left\{ 0 , 1 , \ldots , \varepsilon^{-1} \right\}^d  $, we define $\dis \f_{m}^{\varepsilon}(x)=\prod_{j=1}^{d}f_{m_j}(i_j \varepsilon) \; \;  \text{if}\; \; x\in V_i. $ 
	$ \big\{\, \f_{m}^{\varepsilon}, \, m\in \mathbb{Z}_+^d \,\big\} $ form an orthonormal basis of $ H^{\varepsilon} $ as a subspace of $ L^2\big([0,1]^d\big)$ and are eigenfunctions of $ \Delta_{\varepsilon} $  with eigenvalues\\
	$\dis  -\lambda_{m}^{\varepsilon}=-2\varepsilon^{-2}\sum_{j=1}^{d}\Big( 1 - cos(m_j\pi\varepsilon) \Big). $
	Note that $ \dis \lambda_{m}^{\varepsilon} \longrightarrow \lambda_{m}$ as $\varepsilon \to 0$.\\
	$\bullet$ Basic calculations show that there exists a constant $ c$,  such that for each  $m_j$, $\varepsilon^{-2}\big( 1 - cos(\pi m_j \varepsilon) \big) > c\, m_j^2 .$ \\
	$\bullet$ $ \Delta_{\varepsilon} $ generates a contraction  semigroup $ \mathsf{T}_{\ep}(t) :=\exp \big(\Delta_{\varepsilon} t\big) $  represented  on $ H^{\varepsilon} $ by 
	\begin{eqnarray}\label{rept}
	\mathsf{T}_{\ep}(t)\varphi= \sum_m\exp(-\lambda_{m}^{\varepsilon} t )\langle \; \varphi , \f_{m}^{\varepsilon} \;  \rangle \f_{m} ^{\varepsilon} ,
	\end{eqnarray}
	where the summation is taken on the $\ep^{-d}$ eigenvectors of $\Delta_{\ep}$.
	Note that both $ \Delta_{\varepsilon}$ and $\mathsf{T}_{\ep}(t)$ are self-adjoint and that $ \mathsf{T}_{\ep}(t)\Delta_{\varepsilon}\varphi=\Delta_{\varepsilon}\mathsf{T}_{\ep}(t)\varphi.$ Note also, for any $ J\in\{S, I, R\}$, the semigroup generated by $ \mu_J\Delta$ is $\mt(\mu_Jt)$. In the sequel, we will use the notation $ \mtj(t):=\mt(\mu_Jt)$ and similarly, in the discrete case, we will use the notation  $\mtepj(t):=\mtep(\mu_Jt)$. Also, for any $J\in \{S, I, R\}$, we let $  \lambda_{m,J}:=\mu_J\lambda_{m}$  and $ \lambda_{m,J}^{\ep}:=\mu_J\lambda_{m}^{\varepsilon}. $\\
	$\bullet\;$ We use $  \bV \varphi \bV_{\infty}\! \!:= \un{x\in D}{\sup}\big\vert \varphi(x) \big\vert$ to denote the supremum norm of $\varphi$ in $D$,
	and we define $ \Big\Vert  \bigg( \begin{array}{cr}
	\varphi\\
	\phi
	\end{array}\bigg) \Big\Vert_{\infty} \!\!:=~\bV \varphi \bV_{\infty} +~\bV  \phi \bV_{\infty}. $\\ 
	$\bullet\;$  If $Z$ is a space-time function, we use the notation $ Z(t)=Z(t,.)$. \\
	$\bullet\; $ For $ n \ge 1  $, $ C^n(D)$ denotes the space of real valued continuous functions on $ D$  with  continuous partial derivatives of all orders from 1 to $n$ . We use the standard partial ordering of $\R^d$ and the classical notations:\\
	$u\le v $ if, for all $1\le i\le d$,  $u_i\le v_i$.
	
	\subsubsection{Existence and uniqueness }
	
	Let us set
	$ \dis X_{\ep}= \big( S_{\varepsilon}, I_{\varepsilon}, R_{\varepsilon}
	\big)^T.$  We introduce the function 
	$G : (x;u,v,w) \longmapsto  \left( \begin{array}{cl}
	-\dfrac{\beta(x)\,u\,v}{u+v+w}  \\
	\dfrac{\beta(x)\,u\,v}{u+v+w}-\alpha(x)\, v\\
	\alpha(x)\, v
	\end{array}\right)$.  \\ We use the notation   $ \widetilde{\Delta}_{\varepsilon} X_{\varepsilon} = \Big( 
	\mu_S\Delta_{\varepsilon}S_{\varepsilon},\,
	\mu_I\Delta_{\varepsilon}I_{\varepsilon},\,
	\mu_R\Delta_{\varepsilon}R_{\varepsilon}
	\Big)^T $. Then the compact form of system (\ref{eqdet})  is
	
	\begin{equation}
	\hspace{-3cm}\left \{
	\begin{aligned}\label{fcds1}
	\dfrac{dX_{\varepsilon}}{dt}(t,x_i) & =   \widetilde{\Delta}_{\varepsilon} X_{\varepsilon}(t,x_i) + G\big(x_i ; X_{\varepsilon}(t,x_i)\big), \; \; (t,x_i) \in (0,T)\times D_{\varepsilon} \\
	X_{\varepsilon}(t,x_i)& = X_{\varepsilon}(t,y_i) , \; \text{for} \; x_i\in \partial D_{\varepsilon},\;  x_i \sim y_i \; \text{and} \;  y_i \in \partial_{\vec{n}.out} D_{\varepsilon}\\
	&\hspace{-1.4cm}  X(0,x_i) \ge 0  \; \; \text{and} \; \; 0<S_{\varepsilon}(0,x_i)+I_{\varepsilon}(0,x_i)+R_{\varepsilon}(0,x_i)\le M, \;  \forall \,x_i \in D_\eps.
	\end{aligned}
	\right.
	\end{equation} 
	
	\vspace{-0.2cm}
	
	We have
	
	\begin{lem}\label{ro}
		For each $ \varepsilon >0 $ fixed, the system (\ref{fcds1})  has a unique non-negative solution  $ X_{\varepsilon}\in C^1\big( \mathbb{R}_+ ; \mathbb{R}^{3\varepsilon^{-d}}_+ \big)$. Moreover   $\un{0\le t\le T}{\sup}\Big\Vert X_{\ep}(t)\Big\Vert_{\infty}\le C(\bar{\alpha},\bar{\beta})$.
	\end{lem}
	\begin{pr} 
		Let us define 
		$g :\mathbb{R}^3 \longrightarrow \mathbb{R} $ \\
		
		\vspace{-0.7cm}
		
		\hspace{3.5cm}$(u,v,w)\longmapsto g(u,v,w)=\dfrac{uv}{u+v+w}. $\\ 
		
		We have $  G(x ; u,v,w)= \left( \begin{array}{cl}
		-\beta(x)\, g(u,v,w) \\
		\beta(x)\, g(u,v,w)-\alpha(x)\, v\\
		\alpha(x) \, v
		\end{array}\right).$ We introduce \\
		
		\vspace{0.3cm}
		
		$ g^+(u,v,w)\!=\!\! \left\{ \!\begin{array}{rl}
		\dfrac{u^+v^+}{u^++v^++w^+}, & \mbox{if} \, u^++v^++w^+>0 ;\\
		0, & \mbox{otherwise},
		\end{array}\right. $ 
		where we used the notation $ u^+=\sup(u,0)$. We set
		
		\vspace{0.3cm}
		
		$\dis G^+(x ; u,v,w)= \left( \begin{array}{cl}
		-\beta(x)\, g^+(u,v,w) \\
		\beta(x)\, g^+(u,v,w)-\alpha(x)\, v\\
		\alpha(x) \, v
		\end{array}\right).$
		Let us consider the system of ODEs
		\begin{equation}\label{lc}
		\hspace{-3.5cm}\left \{
		\begin{aligned}
		\dfrac{dX_{\varepsilon}}{dt}(t,x_i) & = \widetilde{\Delta}_{\varepsilon} X_{\varepsilon}(t,x_i) + G^+\big(x_i ; X_{\varepsilon}(t,x_i)\big), \; \; (t,x_i) \in (0,T)\times D_{\varepsilon} \\
		X_{\varepsilon}(t,x_i)& = X_{\varepsilon}(t,y_i) , \; \text{for} \; x_i\in \partial D_{\varepsilon},  x_i \sim y_i \; \text{and} \;  y_i \in \partial_{\vec{n}.out} D_{\varepsilon}\\
		&\hspace{-1.4cm}  X(0,x_i) \ge 0  \; \; \text{and} \; \; 0<S_{\varepsilon}(0,x_i)+I_{\varepsilon}(0,x_i)+R_{\varepsilon}(0,x_i)\le M, \;  \forall \,x_i \in D_\eps.
		\end{aligned}
		\right.
		\end{equation}
		Since $ G^+ $ is globally Lipschitz and $ \widetilde{\Delta}_{\varepsilon}$ is linear, then by the Picard-Lindelöf  theorem the Cauchy problem (\ref{lc}) has a unique solution  $ \tilde{X}_{\varepsilon}\in C^1\big( \mathbb{R}_+ ; \mathbb{R}^{3\varepsilon^{-d}} \big)$. Now let us show that $X_{\varepsilon}(t)\ge0$ for all $t\ge 0.$ Let us set 
		$ t_1=\inf\{t>0 : \text{there exists an index i such that} \; S_{\ep}(t,x_i)< 0\}.$  If $t_1<\infty$ then there exists $i_1$ such that $S_{\ep}(t_1,x_{i_1})=0 $ and for all $j\neq i_1$  $S_{\ep}(t_1,x_j)\ge 0. $ So that $\dfrac{dS_{\varepsilon}}{dt}(t_1,x_{i_1})\ge  0 .$ Thus, either there exists an index $j$ such $x_j\sim x_{i_1}$ and $ S_{\ep}(t_1,x_j)>0$ or  else  $ S_{\ep}(t_1,x_j)=0 $ for all $x_j \sim x_{i_1}$. \\ 
		\textbf{(i)} In the first case  $\dfrac{dS_{\varepsilon}}{dt}(t_1,x_{i_1})> 0$, which contradicts the definition of $ t_1$. \\
		\textbf{(ii)} Let us set $ \dis \mathsf{I}_1=\{\;x_i\in D_{\ep} : S_{\ep}(t_1,x_i)>0\;\}$. If $ \mathsf{I}_1=\varnothing $ then $S_{\ep}(t_1,x_i)=0$, for all $x_i\in D_{\ep} $. The uniqueness of the solution entails that the null vector is solution for the equations satisfied by $S_{\ep}$ on the time interval~$[t_1,+\infty)$. 
		We now assume  that $\mathsf{I}_1\ne\varnothing$, and define
		
		\vspace{-0.5cm}
		
		$$ \mathscr{V}_1(\mathsf{I}_1)=\{\; x_i\in D_{\ep}: x_i\notin \mathsf{I}_1  ,\;  \exists\, x_j \in \mathsf{I}_1 \,\text{such that}\, x_i\sim x_j \;\},$$
		$$\mathscr{V}_2(\mathsf{I}_1)=\{\; x_i\in D_{\eps}: x_i\notin \mathscr{V}_1(\mathsf{I}_1)\cup \mathsf{I}_1,\; \exists\, x_j \in \mathsf{I}_1 \,\text{such that}\, x_i\sim x_j \;\},$$
		
		\vspace{-1cm}
		
		$$\vdots$$
		
		\vspace{-1cm}
		
		$$\mathscr{V}_{k}(\mathsf{I}_1)=\{\; x_i\in D_{\eps}: x_i\notin \mathscr{V}_{k-1}(\mathsf{I}_1)\cup\cdots\cup\mathscr{V}_1(\mathsf{I}_1)\cup \mathsf{I}_1,\;\exists\, x_j \in \mathscr{V}_{k-1}(\mathsf{I}_1) \,\text{such that}\, x_i\sim x_j\;\},\;  k\ge 1.$$
		
		First, note that there exists a positive integer $k$ such that $ \mathscr{V}_{k}(\mathsf{I}_1)= \varnothing$, because there is a finite number of sites. Now, if $ x_i\in \mathscr{V}_1(\mathsf{I}_1)$, then $ S_{\eps}(t_1, x_i)=0$  and $\dfrac{dS_{\eps}}{dt}(t_1, x_i)>0$, which contradicts the definition of $t_1$. Else, assume that $x_i\in \mathscr{V}_2(\mathsf{I}_1)$.  On the one hand,  we have $ S_{\eps}(t_1, x_i)=0$, $\dfrac{dS_{\eps}}{dt}(t_1, x_i)=0$. On the other hand,
		
		\vspace{-0.7cm}
		
		\begin{eqnarray}
		\text{since}\; \; \dfrac{d^2\,S_{\varepsilon}}{dt^2}(t_1,x_i)& = &- \beta(x_i)\dfrac{ I_{\varepsilon}(t_1,x_i)\frac{dS_{\varepsilon}}{dt}(t_1,x_i)+S_{\eps}(t_1,x_i)\frac{dI_{\varepsilon}}{dt}(t_1,x_i)}{S_{\varepsilon}(t_1,x_i)+I_{\varepsilon}(t_1,x_i)+R_{\varepsilon}(t_1,x_i)}\nonumber\\
		&+&\beta(x_i)\dfrac{S_{\eps}(t_1,x_i)I_{\eps}(t_1,x_i)\Big( \frac{dS_{\varepsilon}}{dt}(t_1,x_i)+\frac{dI_{\varepsilon}}{dt}(t_1,x_i)+\frac{dR_{\varepsilon}}{dt}(t_1,x_i) \Big)}{\Big(S_{\varepsilon}(t_1,x_i)+I_{\varepsilon}(t_1,x_i)+R_{\varepsilon}(t_1,x_i)\Big)^2}+ \mu_S\,\Delta_{\varepsilon} \dfrac{dS_{\varepsilon}}{dt}(t_1,x_i), \nonumber
		\end{eqnarray}
		
		\vspace{-0.5cm}
		
		then $ \dis \dfrac{d^2\,S_{\varepsilon}}{dt^2}(t_1,x_i) =  \mu_S\,\Delta_{\varepsilon} \dfrac{dS_{\varepsilon}}{dt}(t_1,x_i)>0$, because $ x_i\in \mathscr{V}_2(\mathsf{I}_1)$, and we obtain again a contradiction. If $x_i\in\mathscr{V}_j(\mathsf{I}_1),$ for $j\ge 2$, we iterate the above argument to obtain  a contradiction.  Then in all cases we obtain a contradiction. So that  $ t_1= \infty.$ Thus $ S_{\eps}(t,x_i)\ge 0$ for  all $ (t,x_i) \in [0,T]\in D_{\eps}.$  Similar arguments hold for  $I_{\varepsilon}$ and $R_{\varepsilon}.$  It follows from the positivity of the solution and the fact that  $ G = G^+ $ on $ \mathbb{R}_+^3$, that  the system (\ref{fcds1}) has a unique global solution $ X_{\varepsilon}\in C^1\big( \mathbb{R}_+ ; \mathbb{R}_+^{3\varepsilon^{-d}} \big)$. Furthermore, writing the solution of (\ref{fcds1}) in its mild semigroup  form, and  using successively the fact that $ \widetilde{\mathsf{T}}_{\eps}:=(
		\mteps,
		\mtepi,
		\mtepr
		)^T$ is a contraction semigroup on $\Big( H^{\varepsilon} , \big(L^{\infty}(D)\big)^3 \Big)$, the fact that $ X(0,x_i)\le 1$ for all $x_i \in D_\eps$ and applying Gronwall's Lemma, we easily obtain that $\un{0\le t\le T}{\sup}\Big\Vert X(t)\Big\Vert_{\infty}\le C(\bar{\alpha},\bar{\beta})$.  This  concludes the proof of the~lemma.
		\fpr
	\end{pr}
	Let us now define  $\; \dis  \mathcal{S}_{\varepsilon}(t,x) = \sum_{i=1}^{\varepsilon^{-d}} S_{\varepsilon}(t,x_i)\mathbf{1}_{V_i}(x),\;$ $\; \dis \mathcal{I}_{\varepsilon}(t,x) = \sum_{i=1}^{\varepsilon^{-d}} I_{\varepsilon}(t,x_i)\mathbf{1}_{V_i}(x), \;\; $
	$\dis \mathcal{R}_{\varepsilon}(t,x) = \sum_{i=1}^{\varepsilon^{-d}} R_{\varepsilon}(t,x_i)\mathbf{1}_{V_i}(x),$   $\dis \beta_{\varepsilon}(t,x) = \sum_{i=1}^{\varepsilon^{-d}}\beta(t,x_i)\mathbf{1}_{V_i}(x),$ \;  $\dis \alpha_{\varepsilon}(t,x) = \sum_{i=1}^{\varepsilon^{-d}} \alpha(t,x_i)\mathbf{1}_{V_i}(x),$   
	and we set $ \dis  \mathcal{X}_{\varepsilon}= \big(
	\mathcal{S}_{\varepsilon},  \mathcal{I}_{\varepsilon},  \mathcal{R}_{\varepsilon}
	\big)^T .$  \\ Note that the previous lemma is true with $\mathcal{X}_{\varepsilon}$  in place of $ X_{\varepsilon}$. Let us set $ \dis  X= \big(\mathbf{s},    \mathbf{i}, \mathbf{r}\big)^T .$
	Then the compact form of the model (\ref{cdm}) is 
	\begin{equation}\label{sytc}
	\hspace{-6cm}\left \{
	\begin{aligned}
	&\hspace{-1.8cm} \dfrac{\partial X}{\partial t}(t,x) =  \widetilde{\Delta}X(t,x) + G\big(x ; X(t,x)\big), \; \; (t,x) \in [0,T]\times D \\
	\dfrac{\partial X}{\partial n_{\text{out}}}(t,x) &  = 0,  \;  \; \; \text{for}\; x\in \partial D\\
	&\hspace{-1.8cm}  X(0,x)\ge 0 \; \text{and} \; 0< \mathbf{s}(0,x)+\mathbf{i}(0,x)+\mathbf{r}(0,x)\le M.
	\end{aligned}
	\right.
	\end{equation}
	
	\vspace{-0.5cm}
	
	Writing the solution of (\ref{sytc}) in its mild semigroup  form, we have $\dis X(t) = \widetilde{\mt}(t)X(0)+ \int_0^t\widetilde{\mt}(t-r)G\big(X(r)\big)dr$, where we used the notation
	$ \dis \widetilde{\mt}(t)X:= \left( \begin{array}{cl}
	\mts(t)\mathbf{s} \\
	\mti(t)\mathbf{i} \\
	\mtr(t)\mathbf{r}
	\end{array}\right)$ and  similarly  for $\widetilde{\mt}(t-r)G\big(X(r)\big).$
	\begin{lem}
		The initial value probem (\ref{sytc}) has a unique  solution $ X\in C\Big([0,T]\,; \big(L^{\infty}(D)\big)^3\Big)$.
		
	\end{lem}
	
	\begin{pr}
		For $0\le u(0,.)\le 1$, we define a mapping F $:\dis C\Big([0,T]\,; \big(L^{\infty}(D)\big)^3\Big)\longrightarrow  C\Big([0,T]\,; \big(L^{\infty}(D)\big)^3\Big)$ by
		\begin{eqnarray}\label{f}
		(Fu)(t) = \widetilde{\mt}(t)u(0)+ \int_0^t\widetilde{\mt}(t-r)G\big(u(r)\big)dr.
		\end{eqnarray}
		
		Let $ u,v  \in C\Big([0,T]\,; \big(L^{\infty}(D)\big)^3\Big)$ such that $u(0)=v(0)$.
		Using the fact that  $\widetilde{\mt}$ is a contraction semigroup on $\big(L^{\infty}(D)\big)^3$ and G is globally Lipschitz, it follows that
		\begin{eqnarray}
		\Big\Vert (Fu)(t)-(Fv)(t)\Big\Vert_{\infty}\le C\int_0^t \Big\Vert u(r)-v(r)\Big\Vert_{\infty}dr,\; \; \text{for all} \; t\in[0,T],\nonumber
		\end{eqnarray}
		where $C$ is the Lipschitz constant of $G$. Then we have
		\begin{eqnarray}\label{if}
		\Big\Vert (Fu)(t)-(Fv)(t)\Big\Vert_{\infty}\le CT \un{0\le r\le T}{\sup} \Big\Vert u(r)-v(r)\Big\Vert_{\infty}. 
		\end{eqnarray}
		Using (\ref{f}), (\ref{if}) and induction on $n$, it follows that 
		\begin{eqnarray}
		\un{0\le t\le T}{\sup}\Big\Vert (F^nu)(t)-(F^nv)(t)\Big\Vert_{\infty}\le \dfrac{(CT)^n}{ n!} \un{0\le t\le T}{\sup} \Big\Vert u(t)-v(t)\Big\Vert_{\infty}. 
		\end{eqnarray}
		For $n$ large enough $F^n$ is a contraction (since $ \dfrac{(CT)^n}{ n!}<1$). It follows from the Banach contraction principle that F has a unique fixed point in $ X\in C\Big([0,T]\,; \big(L^{\infty}(D)\big)^3\Big).$ This fixed point is the solution of 
		\begin{eqnarray}\label{msol}
		X(t) = \widetilde{\mt}(t)X(0)+ \int_0^t\widetilde{\mt}(t-r)G\big(X(r)\big)dr.
		\end{eqnarray}
		The mild solution of (\ref{msol}) is in fact a classical solution of (\ref{sytc}), see Theorem 3.1,  Chapter 7 of Smith~(1995). Note  that the Corollary 3.1,  Chapter 7 of the above reference ensures that $ X(t)\ge 0$,  $\forall \,t  \ge 0$. 
		\fpr
	\end{pr}
	
	\subsubsection{Relation between the system of PDEs and its discretisation}
	
	We will now prove that  $ \mathcal{X}_{\varepsilon}$  converges to X  as the mesh  size $ \varepsilon$  of the grid tends to zero.
	
	\begin{thm}\label{ct}  Let  us consider an initial condition  $ X(0) \in \big(L^{\infty}(D)\big)^3$. For all  $T > 0 $,  \\
		$\dis \un{t\in [0,T]}{\sup}\Big\Vert \mathcal{X}_{\varepsilon}(t)- X(t)\Big\Vert_{\infty} \longrightarrow 0 $, as \,  $ \varepsilon \to 0 . $ 
	\end{thm}
	\begin{pr}
		Using the variation of constants formula, we have 
		\begin{eqnarray}\mathcal{S}_{\varepsilon}(t) 
		&=& \mteps(t)\mathcal{S}_{\varepsilon}(0) -  \int_0^t  \mteps(t-r)\Big[\dfrac{\beta_{\ep}(.)\mathcal{S}_{\varepsilon}(r) \mathcal{I}_{\varepsilon}(r)}{\mathcal{S}_{\varepsilon}(r) +\mathcal{I}_{\varepsilon}(r) +\mathcal{R}_{\varepsilon}(r) }\Big]dr,  \nonumber 
		\end{eqnarray}
		$$  \mathbf{s}(t) = \mts(t)\mathbf{s}(0) -  \int_0^t\mts(t-r)\Big[\dfrac{\beta(.)\mathbf{s}(r) \mathbf{i}(r)}{\mathbf{s}(r)+\mathbf{i}(r)+\mathbf{r}(r)}\Big]dr. $$ 
		Recall that $ \mathcal{S}_{\varepsilon}(0)= P_{\varepsilon}\,\mathbf{s}(0) $, so that
		\begin{eqnarray}
		\mathcal{S}_{\varepsilon}(t) -\mathbf{s}(t)
		&=& \mteps(t) P_{\varepsilon}\,\mathbf{s}(0) -\mts(t)\mathbf{s}(0) -  \int_0^t  \mteps(t-r)\Big[\dfrac{\beta_{\ep}(.)\mathcal{S}_{\varepsilon}(r) \mathcal{I}_{\varepsilon}(r)}{\mathcal{S}_{\varepsilon}(r) +\mathcal{I}_{\varepsilon}(r) +\mathcal{R}_{\varepsilon}(r) }\Big]dr \nonumber\\
		&&+ \int_0^t\mts(t-r)\Big[\dfrac{\beta(.)\mathbf{s}(r) \mathbf{i}(r)}{\mathbf{s}(r)+\mathbf{i}(r)+\mathbf{r}(r)}\Big]dr. \nonumber 
		\end{eqnarray}
		We have
		\begin{eqnarray}\label{eqm}
		\Big\Vert \mathcal{S}_{\varepsilon}(t) - \mathbf{s}(t) \Big\Vert_{\infty}
		&\le& \Big\Vert \mteps(t) P_{\varepsilon}\,\mathbf{s}(0) -\mts(t)\mathbf{s}(0) \Big\Vert_{\infty}\nonumber\\
		&&\hspace{-2cm}+\int_0^t \Big\Vert \mteps(t-r)\Big[\dfrac{\beta_{\ep}(.)\mathcal{S}_{\varepsilon}(r) \mathcal{I}_{\varepsilon}(r)}{\mathcal{S}_{\varepsilon}(r) +\mathcal{I}_{\varepsilon}(r) +\mathcal{R}_{\varepsilon}(r) }\Big]- \mts(t-r)\Big[\dfrac{\beta(.)\mathbf{s}(r) \mathbf{i}(r)}{\mathbf{s}(r)+\mathbf{i}(r)+\mathbf{r}(r)}\Big] \Big\Vert_{\infty} dr \nonumber\\
		&&\hspace{-2cm} \le \Big\Vert \mteps(t) P_{\varepsilon}\,\mathbf{s}(0) -\mts(t)\mathbf{s}(0) \Big\Vert_{\infty}\nonumber\\
		&&\hspace{-2cm} + \int_0^t\Big\Vert \mteps(t-r)\Big[\dfrac{\beta{\ep}(.)\mathcal{S}_{\varepsilon}(r) \mathcal{I}_{\varepsilon}(r)}{\mathcal{S}_{\varepsilon}(r) +\mathcal{I}_{\varepsilon}(r) +\mathcal{R}_{\varepsilon}(r) }\Big]-\mteps(t-r)P_{\varepsilon}\Big[\dfrac{\beta(.)\mathbf{s}(r) \mathbf{i}(r)}{\mathbf{s}(r)+\mathbf{i}(r)+\mathbf{r}(r)}\Big] \Big\Vert_{\infty} dr\nonumber \\
		&& \hspace{-2cm} +\int_0^t \Big\Vert \mteps(t-r)P_{\varepsilon}\Big[\dfrac{\beta(.)\mathbf{s}(r) \mathbf{i}(r)}{\mathbf{s}(r)+\mathbf{i}(r)+\mathbf{r}(r)}\Big]-\mts(t-r)\Big[\dfrac{\beta(.)\mathbf{s}(r) \mathbf{i}(r)}{\mathbf{s}(r)+\mathbf{i}(r)+\mathbf{r}(r)}\Big] \Big\Vert_{\infty} dr.
		\end{eqnarray}
		
		Let us estimate each term of the right-hand side of this inequality. 
		
		Since $ \mathbf{s}(0)\in L^{\infty}(D)$, it then follows from  Kato ( \cite{Kato66} pp. 512-513 ), that
		
		$~\hspace{2cm}\Big\Vert \mteps(t)P_{\varepsilon}\mathbf{s}(0) -\mts(t)\mathbf{s}(0) \Big\Vert_{\infty} \longrightarrow 0, \; \;  \text{uniformly on} \; [0,T].$
		
		Using the fact that $ \mteps$ is a contraction semigroup on $\Big( H^{\varepsilon} , \Vert . \Vert_{\infty}\Big)$, we obtain
		\begin{eqnarray}
		\int_0^t\bigg\Vert \mteps(t-r)\Big[\dfrac{\beta_{\ep}(.)\mathcal{S}_{\varepsilon}(r) \mathcal{I}_{\varepsilon}(r)}{\mathcal{S}_{\varepsilon}(r) +\mathcal{I}_{\varepsilon}(r) +\mathcal{R}_{\varepsilon}(r) }\Big]-\mteps(t-r)P_{\varepsilon}\Big[\dfrac{\beta(.)\mathbf{s}(r) \mathbf{i}(r)}{\mathbf{s}(r)+\mathbf{i}(r)+\mathbf{r}(r)}\Big] \bigg\Vert_{\infty}dr 
		&&\nonumber \\
		&&\hspace{-8cm}\le \bar{\beta}\int_0^t\bigg\Vert \dfrac{\mathcal{S}_{\varepsilon}(r) \mathcal{I}_{\varepsilon}(r)}{\mathcal{S}_{\varepsilon}(r) +\mathcal{I}_{\varepsilon}(r) +\mathcal{R}_{\varepsilon}(r) }-\dfrac{\mathbf{s}(r) \mathbf{i}(r)}{\mathbf{s}(r)+\mathbf{i}(r)+\mathbf{r}(r)} \bigg\Vert_{\infty} dr \nonumber\\
		&&\hspace{-8cm}+\int_0^t  \bigg\Vert P_{\varepsilon}\Big[\dfrac{\beta(.)\mathbf{s}(r) \mathbf{i}(r)}{\mathbf{s}(r)+\mathbf{i}(r)+\mathbf{r}(r)}\Big]- \dfrac{\beta(.)\mathbf{s}(r) \mathbf{i}(r)}{\mathbf{s}(r)+\mathbf{i}(r)+\mathbf{r}(r)} \bigg\Vert_{\infty} dr \nonumber \\
		&&\hspace{-8cm} \le  \bar{\beta}\int_0^t \Bigg( 2 \bigg\Vert \mathcal{S}_{\varepsilon}(r) - \mathbf{s}(r) \bigg\Vert_{\infty} + 2\bigg\Vert \mathcal{I}_{\varepsilon}(r)-\mathbf{i}(r)\bigg\Vert_{\infty}  +\bigg\Vert \mathcal{R}_{\varepsilon}(r) - \mathbf{r}(r) \bigg\Vert_{\infty} \Bigg) dr \nonumber\\
		&&\hspace{-8cm}+ \int_0^t\bigg\Vert P_{\varepsilon}\Big[\dfrac{\beta(.)\mathbf{s}(r) \mathbf{i}(r)}{\mathbf{s}(r)+\mathbf{i}(r)+\mathbf{r}(r)}\Big]- \dfrac{\beta(.)\mathbf{s}(r) \mathbf{i}(r)}{\mathbf{s}(r)+\mathbf{i}(r)+\mathbf{r}(r)} \bigg\Vert_{\infty}  dr. \nonumber
		\end{eqnarray}
		Since $ \dfrac{\beta(.)\mathbf{s}(r) \mathbf{i}(r)}{\mathbf{s}(r)+\mathbf{i}(r)+\mathbf{r}(r)}\in L^{\infty}(D)$, then for the last term  (\ref{eqm}), we have
		$$\bigg\Vert \mteps(t-r)P_{\varepsilon}\Big[\dfrac{\beta(.)\mathbf{s}(r) \mathbf{i}(r)}{\mathbf{s}(r)+\mathbf{i}(r)+\mathbf{r}(r)}\Big]-\mts(t-r)\Big[\dfrac{\beta(.)\mathbf{s}(r) \mathbf{i}(r)}{\mathbf{s}(r)+\mathbf{i}(r)+\mathbf{r}(r)}\Big] \bigg\Vert_{\infty} \longrightarrow 0, \;  \text{uniformly on}\; [0,T],$$ ( \text{thanks by Kato \cite{Kato66}}, \text{chapter 9, Section 3} ).
		Consequently
		$$\hspace{-5cm} \int_0^t \bigg\Vert \mteps(t-r)P_{\varepsilon}\Big[\dfrac{\beta(.)\mathbf{s}(r) \mathbf{i}(r)}{\mathbf{s}(r)+\mathbf{i}(r)+\mathbf{r}(r)}\Big]-\mts(t-r)\Big[\dfrac{\beta(.)\mathbf{s}(r) \mathbf{i}(r)}{\mathbf{s}(r)+\mathbf{i}(r)+\mathbf{r}(r)}\Big] \bigg\Vert_{\infty}dr \longrightarrow 0 .$$
		\begin{eqnarray}
		\Big\Vert \mathcal{S}_{\varepsilon}(t)-\mathbf{s}(t)\Big\Vert_{\infty}
		&\le& a_{\ep}(t)+C(\bar{\beta}) \int_0^t \Bigg(\Big\Vert \mathcal{S}_{\varepsilon}(r) -\mathbf{s}(r) \Big\Vert_{\infty} + \Big\Vert \mathcal{I}_{\varepsilon}(r)-\mathbf{i}(r)+ \Big\Vert \mathcal{R}_{\varepsilon}(r) -\mathbf{r}(r) \Big\Vert_{\infty} \Bigg) dr, \; \; \text{where} \nonumber 
		\end{eqnarray}
		
		\vspace{-0.8cm}	
		
		\begin{eqnarray}
		a_{\ep}(t)&=&\Big\Vert \mteps(t) P_{\varepsilon}\,\mathbf{s}(0) -\mts(t)\mathbf{s}(0) \Big\Vert_{\infty}+ \int_0^t  \bigg\Vert P_{\varepsilon}\Big[\dfrac{\beta(.)\mathbf{s}(r) \mathbf{i}(r)}{\mathbf{s}(r)+\mathbf{i}(r)+\mathbf{r}(r)}\Big]- \dfrac{\beta(.)\mathbf{s}(r) \mathbf{i}(r)}{\mathbf{s}(r)+\mathbf{i}(r)+\mathbf{r}(r)} \bigg\Vert_{\infty} dr \nonumber\\
		&&+\int_0^t \bigg\Vert \mteps(t-r)P_{\varepsilon}\Big[\dfrac{\beta(.)\mathbf{s}(r) \mathbf{i}(r)}{\mathbf{s}(r)+\mathbf{i}(r)+\mathbf{r}(r)}\Big]-\mts(t-r)\Big[\dfrac{\beta(.)\mathbf{s}(r) \mathbf{i}(r)}{\mathbf{s}(r)+\mathbf{i}(r)+\mathbf{r}(r)}\Big] \bigg\Vert_{\infty} dr, \nonumber
		\end{eqnarray}
		
		Exactly in the same way we have a similar inequality for $\Big\Vert \mathcal{I}_{\varepsilon}(t)-\mathbf{i}(t)\Big\Vert_{\infty}$ and $\Big\Vert \mathcal{R}_{\varepsilon}(t)-\mathbf{r}(t)\Big\Vert_{\infty}$ with $\mti$, $\mtr$ in place of $\mts$, and $\mtepi$ , $\mtepr$ in place of $\mteps$, respectively. Combining those estimates we obtain
		\begin{eqnarray}
		\Big\Vert \mathcal{X}_{\varepsilon}(t)-X(t)\Big\Vert_{\infty} 
		&\le & \tilde{a}_{\ep}(t) + C(\bar{\alpha},\bar{\beta}) \int_0^t \Big\Vert \mathcal{X}_{\varepsilon}(r) - X(r) \Big\Vert_{\infty}dr, \nonumber
		\end{eqnarray}
		where $\dis \un{0\le t\le T}{\sup} \tilde{a}_{\ep}(t) \longrightarrow 0,$ as $\ep \to 0.$  Applying Gronwall's Lemma, it follows that 
		
		\vspace{-0.2cm}	
		
		\begin{eqnarray}\label{nx}
		\un{0\le t\le T}{\text{sup}}\Big\Vert \mathcal{X}_{\varepsilon}(t)-X(t)\Big\Vert_{\infty}
		&\le & \un{0\le t\le T}{\sup} \tilde{a}_{\ep}(t) e^{C(\bar{\alpha},\bar{\beta})T}. \nonumber
		\end{eqnarray}
		
		\vspace{-0.2cm}	
		
		Finally, the theorem follows from the fact that the right-hand side  tends to zero as $ \varepsilon \to 0 $. 
		\fpr
	\end{pr}
	\subsection{The stochastic model}
	\rhead{The stochastic model}
	
	Deterministic models describe the spread of disease under the assumptions of mass action, relying on the law of large numbers. The most natural way  to describe the spread of disease is stochastic.
	The previous models are based on the hypothesis of a population of large size. When it is not the case, the interactions between the individuals are not uniform but possess an intrinsic random character. We are going to expose now a probabilistic version of the previous model.  
	For each given site, Poisson processes  count the number of new infections, removal and  migrations between sites during  time. So the propagation of the illness can be modeled by  the following  system of stochastic differential equations
	\begin{equation}
	\hspace{-0.3cm}\left\{ 
	\begin{aligned} 
	\bigskip 
	S^{\varepsilon}(t,x_i) &=  S^{\varepsilon}(0,x_i) - \mathrm{P}_{x_i}^{inf}\left(  \int_0^t \dfrac{\beta(x_i)S^{\varepsilon}(r,x_i)I^{\varepsilon}(r,x_i)}{S^{\varepsilon}(r,x_i)+I^{\varepsilon}(r,x_i)+R^{\varepsilon}(r,x_i)}dr \right) \\
	& -  \sum_{y_i\sim x_i}\mathrm{P}_{S,x_i,y_i}^{mig}\left( \int_0^t \frac{\mu_S }{\varepsilon^2}S^{\varepsilon}(r,x_i)dr \right) + \sum_{y_i\sim x_i}\mathrm{P}_{S,y_i,x_i}^{mig}\left( \int_0^t \frac{\mu_S}{\varepsilon^2} S^{\varepsilon}(r,y_i)dr \right)
	\\[2mm]
	I^{\varepsilon}(t,x_i) & =  I^{\varepsilon}(0,x_i) +\mathrm{P}_{x_i}^{inf}\left(\int_0^t \dfrac{\beta(x_i)S^{\varepsilon}(r,x_i)I^{\varepsilon}(r,x_i)}{S^{\varepsilon}(r,x_i)+I^{\varepsilon}(r,x_i)+R^{\varepsilon}(r,x_i)}dr \right)  - \mathrm{P}_{x_i}^{rec}\left(   \int_0^t \alpha(x_i)I^{\varepsilon}(r,x_i)dr \right) \\
	&  - \sum_{y_i\sim x_i}\mathrm{P}_{I,x_i,y_i}^{mig}\left( \int_0^t \frac{\mu_I}{\varepsilon^2} I^{\varepsilon}(r,x_i)dr \right) + \sum_{y_i\sim x_i}\mathrm{P}_{I,y_i,x_i}^{mig}\left( \int_0^t\frac{\mu_I}{\varepsilon^2} I^{\varepsilon}(r,y_i)dr \right)
	\\[2mm] 
	R^{\varepsilon}(t,x_i) & =  R^{\varepsilon}(0,x_i) + \mathrm{P}_{x_i}^{rec}\left(   \int_0^t \alpha(x_i)I^{\varepsilon}(r,x_i)dr \right)-\sum_{y_i\sim x_i}\mathrm{P}_{R,x_i,y_i}^{mig}\left( \int_0^t \frac{\mu_R}{\varepsilon^2} R^{\varepsilon}(r,x_i)dr \right) \\
	& + \sum_{y_i\sim x_i}\mathrm{P}_{R,y_i,x_i}^{mig}\left( \int_0^t \frac{\mu_R}{\varepsilon^2} R^{\varepsilon}(r,y_i)dr \right), \quad (t,x_i) \in [0,T]\times D_{\varepsilon} ,
	\end{aligned}
	\right. 
	\end{equation}
	where all the $\mathrm{P}_j$'s are mutually independent standard Poisson processes.
	In this system, at a given site $ x_i$ 
	\begin{enumerate}
		\item[$\bullet $] infection of a susceptible happens at rate $\dis \beta(x_i)\dfrac{S^{\varepsilon}(t,x_i)}{S^{\varepsilon}(t,x_i)+I^{\varepsilon}(t,x_i)+R^{\varepsilon}(t,x_i)}I^{\varepsilon}(t,x_i)$. Then\\ $ \dis  \mathrm{P}_{x_i}^{inf}\left(\int_0^t \dfrac{\beta(x_i)S^{\varepsilon}(r,x_i)I^{\varepsilon}(r,x_i)}{S^{\varepsilon}(r,x_i)+I^{\varepsilon}(r,x_i)+R^{\varepsilon}(r,x_i)}dr \right) $ counts the number of transitions of type $ S^{\varepsilon} \longrightarrow I^{\varepsilon} $ at site $x_i$ between time $0$ and time $ t$ .
		\item[$\bullet $] recovery of an infectious happens at rate $ \alpha(x_i) I^{\varepsilon}(t,x_i) $, so $ \dis \mathrm{P}_{x_i}^{rec}\left(\int_0^t \alpha(x_i) I^{\varepsilon}(r,x_i)dr \right)$ counts the number of  transitions of type  $ I^{\varepsilon} \longrightarrow R^{\varepsilon} $ at site $x_i$ between time $0$ and time $ t$.  
		\item[$\bullet $] The term $ \dis \mathrm{P}_{S,x_i,y_i}^{mig}\left( \int_0^t \frac{\mu_S }{\varepsilon^2}S^{\varepsilon}(r,x_i)dr \right) $ counts  the number of migrations of susceptibles from site $x_i$ to  $y_i$ (where  $x_i$ and $y_i$ are neighbors), those events happen at rate $ \dfrac{\mu_S}{\varepsilon^2} S^{\varepsilon}(t,x_i)$ ; and similarly for the compartments $I^{\varepsilon} $ and $ R^{\varepsilon}.$  
	\end{enumerate}
	We introduce the martingales $ \mathrm{M}_j(t) = \mathrm{P}_j(t) - t $ and we look instead at the renormalized model by dividing the number of individuals in each compartment and at each site by $ \mathbf{N} $. Hence by setting 
	\begin{eqnarray}\label{rn}
	S_{\mathbf{N},\varepsilon}(t,x_i) = \dfrac{S^{\varepsilon}(t,x_i)}{\mathbf{N}},\; \; I_{\mathbf{N},\varepsilon}(t,x_i) = \dfrac{I^{\varepsilon}(t,x_i)}{\mathbf{N}}, \; \; \text{and} \; \; R_{\mathbf{N},\varepsilon}(t,x_i) = \dfrac{R^{\varepsilon}(t,x_i)}{\mathbf{N}},
	\end{eqnarray} the equations in the various compartments read
	\begin{equation}\label{eqdifs}
	\hspace{-0.8cm}\left\{ 
	\begin{aligned} 
	\bigskip
	S_{\mathbf{N},\varepsilon}(t,x_i)
	& =  S_{\mathbf{N},\varepsilon}(0,x_i) - \int_0^t \dfrac{\beta(x_i) S_{\mathbf{N},\varepsilon}(r,x_i)I_{\mathbf{N},\varepsilon}(r,x_i)}{S_{\mathbf{N},\varepsilon}(r,x_i)+I_{\mathbf{N},\varepsilon}(r,x_i)+R_{\mathbf{N},\varepsilon}(r,x_i)}dr\\
	&\hspace{-1.5cm}+ \int_0^t \mu_S\Delta_{\varepsilon} S_{\mathbf{\mathbf{N}},\varepsilon}(r,x_i)dr 
	-\frac{1}{\mathbf{N}}\mathrm{M}_{x_i}^{inf}\left(\mathbf{N}\int_0^t \dfrac{\beta(x_i) S_{\mathbf{N},\varepsilon}(r,x_i)I_{\mathbf{N},\varepsilon}(r,x_i)}{S_{\mathbf{N},\varepsilon}(r,x_i)+I_{\mathbf{N},\varepsilon}(r,x_i)+R_{\mathbf{N},\varepsilon}(r,x_i)}dr \right)\\
	&\hspace{-1.5cm}-\sum_{y_i\sim x_i}\frac{1}{\mathbf{N}}\mathrm{M}_{S,x_i,y_i}^{mig}\left( \mathbf{N}\int_0^t \dfrac{\mu_S }{\varepsilon^2}  S_{\mathbf{N},\varepsilon}(r,x_i)dr \right) 
	+ \sum_{y_i\sim x_i}\frac{1}{\mathbf{N}}\mathrm{M}_{S,y_i,x_i}^{mig}\left(  \mathbf{N}\int_0^t \dfrac{\mu_S }{\varepsilon^2}S_{\mathbf{N}, \varepsilon}(r,y_i)dr \right)  
	\\[3mm]
	I_{\mathbf{N},\varepsilon}(t,x_i)
	& = I_{\mathbf{N},\varepsilon}(0,x_i) + \int_0^t \dfrac{\beta(x_i) S_{\mathbf{N},\varepsilon}(r,x_i)I_{\mathbf{N},\varepsilon}(r,x_i)}{S_{\mathbf{N},\varepsilon}(r,x_i)+I_{\mathbf{N},\varepsilon}(r,x_i)+R_{\mathbf{N},\varepsilon}(r,x_i)}dr -  \int_0^t  \alpha(x_i) \, I_{\mathbf{N},\varepsilon}(r,x_i)dr\\
	&\hspace{-1.5cm}+\mu_I\int_0^t \Delta_{\varepsilon} I_{\mathbf{N},\varepsilon}(r,x_i)dr +\frac{1}{\mathbf{N}}\mathrm{M}_{x_i}^{inf}\left(\mathbf{N}\int_0^t\dfrac{\beta(x_i) S_{\mathbf{N},\varepsilon}(r,x_i)I_{\mathbf{N},\varepsilon}(r,x_i)}{S_{\mathbf{N},\varepsilon}(r,x_i)+I_{\mathbf{N},\varepsilon}(r,x_i)+R_{\mathbf{N},\varepsilon}(r,x_i)}dr \right) \\
	&\hspace{-1.5cm}- \frac{1}{\mathbf{N}}\mathrm{M}_{x_i}^{rec}\left(\mathbf{N}\int_0^t \alpha (x_i)I_{\mathbf{N},\varepsilon}(r,x_i)dr \right) 
	-\sum_{y_i\sim x_i}\frac{1}{\mathbf{N}}\mathrm{M}_{I,x_i,y_i}^{mig}\left(\mathbf{N} \int_0^t \dfrac{\mu_I }{\varepsilon^2}I_{\mathbf{N},\varepsilon}(r,x_i)dr \right) \\
	&\hspace{-1.5cm}+ \sum_{y_i\sim x_i}\frac{1}{\mathbf{N}}\mathrm{M}_{I,y_i,x_i}^{mig}\left(\mathbf{N}\int_0^t \dfrac{\mu_I}{\varepsilon^2}I_{\mathbf{N},\varepsilon}(r,y_i)dr \right)
	\\[3mm]
	R_{\mathbf{N},\varepsilon}(t,x_i)
	& = R_{\mathbf{N},\varepsilon}(0,x_i) +  \int_0^t  \alpha(x_i) \, I_{\mathbf{N},\varepsilon}(r,x_i)dr +\int_0^t \mu_R \Delta_{\varepsilon} R_{\mathbf{N},\varepsilon}(r,x_i)dr \\
	&+ \frac{1}{\mathbf{N}}\mathrm{M}_{x_i}^{rec}\left( \mathbf{N}\int_0^t \alpha(x_i)I_{\mathbf{N},\varepsilon}(r,x_i)dr \right) 
	-\sum_{y_i\sim x_i}\frac{1}{\mathbf{N}}\mathrm{M}_{R,x_i,y_i}^{mig}\left( \mathbf{N}\int_0^t \dfrac{\mu_R}{\varepsilon^2} R_{\mathbf{N},\varepsilon}(r,x_i)dr \right)\\
	&+ \sum_{y_i\sim x_i}\frac{1}{\mathbf{N}}\mathrm{M}_{R,y_i,x_i}^{mig}\left(\mathbf{N}\int_0^t \dfrac{\mu_R}{\varepsilon^2}R_{\mathbf{N},\varepsilon}(r,y_i)dr \right). 
	\end{aligned}
	\right. 
	\end{equation}
	
	Let  $S_{\mathbf{N},\varepsilon}(t)$ and $I_{\mathbf{N},\varepsilon}(t)$ and $R_{\mathbf{N},\varepsilon}(t)$ denote respectively the vectors which describe the "proportions" of susceptibles, infectious and removed in the population at the various sites at time t :\\
	
	$ S_{\mathbf{N},\varepsilon}(t) = \left( \begin{array}{cl}
	S_{\mathbf{N},\varepsilon}(t,x_1)\\
	\vdots \\
	S_{\mathbf{N},\varepsilon}(t,x_{\ell})
	\end{array}
	\right)$,  \quad
	$ I_{\mathbf{N},\varepsilon}(t) = \left( \begin{array}{cl}
	I_{\mathbf{N},\varepsilon}(t,x_1)\\
	\vdots \\
	I_{\mathbf{N},\varepsilon}(t,x_{\ell}) 
	\end{array} 
	\right) $  and \; \; 
	$ R_{\mathbf{N},\varepsilon}(t) = \left( \begin{array}{cl}
	R_{\mathbf{N},\varepsilon}(t,x_1)\\
	\vdots \\
	R_{\mathbf{N},\varepsilon}(t,x_{\ell})
	\end{array}
	\right),$
	
	where $ \ell $ is the total number of locations. Let us set
	$\dis Z_{\mathbf{N},\varepsilon}(t)= \left( \begin{array}{cl}
	S_{\mathbf{N},\varepsilon}(t)\\
	I_{\mathbf{N},\varepsilon}(t)\\
	R_{\mathbf{N},\varepsilon}(t)
	\end{array}
	\right);  $
	then the aggregated form of the stochastic model is
	
	\vspace{-01cm}
	
	\begin{eqnarray}\label{fag}
	\hspace{-1cm} Z_{\mathbf{N},\varepsilon}(t) 
	&=& Z_{\mathbf{N},\varepsilon}(0)+ \int_0^t b_{\varepsilon}\Big(r,Z_{\mathbf{N},\varepsilon}(r)\Big)dr + \sum_{j=1}^{k_\eps}\frac{h_j}{\mathbf{N}} M_j\Bigg( \mathbf{N}\int_0^t \beta_j\big(r, Z_{\mathbf{N},\varepsilon}(r)\big)dr\Bigg),
	\end{eqnarray}
	
	\vspace{-0.6cm}
	
	where $ \dis  \forall\, r \,\geq \,0 , \;  b_{\varepsilon}\Big(r,Z_{\mathbf{N},\varepsilon}(r)\Big)\!=\! \!  \sum_{j=1}^{k_\ep}\!h_j \beta_j\Big(r,Z_{\mathbf{N},\varepsilon}(r)\Big) \; ; $ the coordinates of each vector $h_j$ are either $ -1$, $0$ or~$1$ and 
	
	\vspace{-1cm}
	
	\begin{eqnarray} 
	\beta_j\Big(r,Z_{\mathbf{N},\varepsilon}(r)\Big)&\!\!\!\!\in\!\!\!\!&\left\{ \dfrac{\beta(.) S_{\mathbf{N},\varepsilon}(r,x_i)I_{\mathbf{N},\varepsilon}(r,.)}{S_{\mathbf{N},\varepsilon}(r,.)+I_{\mathbf{N},\varepsilon}(r,.)+R_{\mathbf{N},\varepsilon}(r,.)}, \; \dfrac{\mu_S}{\varepsilon^2}S_{\mathbf{N},\varepsilon}(r,.) , \;   \dfrac{\mu_I}{\varepsilon^2} I_{\mathbf{N},\varepsilon}(r,.) , \;  \dfrac{\mu_R}{\varepsilon^2} R_{\mathbf{N},\varepsilon}(r,.) , \;   \alpha(.) I_{\mathbf{N},\varepsilon}(r,.)\right\}, \nonumber
	\end{eqnarray}
	$k_\ep$ is the total number of Poisson processes in the model. Note that $ \dis  b_{\varepsilon}\big( r, Z_{\mathbf{N},\varepsilon}(r)\big)= \widetilde{\Delta}_{\ep}Z_{\mathbf{N},\varepsilon}(r)+ G\big(Z_{\mathbf{N},\varepsilon}(r)\big), $ where \\
	$ \widetilde{\Delta}_\ep Z_{\mathbf{N},\varepsilon}(r)= \left(\begin{array}{lc} 
	\mu_S\Delta_{\varepsilon}S_{\mathbf{N},\varepsilon}(r) \\
	\mu_I\Delta_{\varepsilon}I_{\mathbf{N},\varepsilon}(r)\\
	\mu_R\Delta_{\varepsilon}R_{\mathbf{N},\varepsilon}(r)
	\end{array}  \right)  $ \;  and \; 
	$ G\big(Z_{\mathbf{N},\varepsilon}(r)\big)= \left(\begin{array}{lc} 
	- \dfrac{\beta(.) S_{\mathbf{N},\varepsilon}(r)I_{\mathbf{N},\varepsilon}(r)}{S_{\mathbf{N},\varepsilon}(r)+I_{\mathbf{N},\varepsilon}(r,x_i)+I_{\mathbf{N},\varepsilon}(r)} \\
	\\
	\bigskip
	
	\dfrac{\beta(.) S_{\mathbf{N},\varepsilon}(r)I_{\mathbf{N},\varepsilon}(r)}{S_{\mathbf{N},\varepsilon}(r)+I_{\mathbf{N},\varepsilon}(r)+I_{\mathbf{N},\varepsilon}(r)} - \alpha(.) I_{\mathbf{N},\varepsilon}(r)\\
	\bigskip
	\alpha(.) I_{\mathbf{N},\varepsilon}(r)
\end{array}  \right).  $\\

\begin{Large}\textbf{Existence and uniqueness}\end{Large}\\ \\
At the begining of the epidemic,  the proportions of the population in various compartments take their values in the  discrete set $\big\{\; \dfrac{n}{\mathbf{N}}, \;  n=0,1,\cdots \big\}$, and since the Poisson processes are mutually independent, this implies that the components of $Z_{\mathbf{N},\varepsilon}(t)$ remain non-negative for all $ t\ge 0$. Indeed, let us consider for example the component $ S_{\mathbf{N},\varepsilon}$. Since all jumps of each $S_{\mathbf{N},\varepsilon}(t,x_i)$ are of size $\pm \dfrac{1}{\mathbf{N}}$, before becoming negative, $S_{\mathbf{N},\varepsilon}(t,x_i)$ is zero. But as long as $S_{\mathbf{N},\varepsilon}(t,x_i)=0$, the rate of its negative jumps is zero, hence $S_{\mathbf{N},\varepsilon}(t,x_i)$ cannot become negative.
\\
$\dis \sum_{i=1}^{\varepsilon^{-d}}\Big( S_{\mathbf{N},\varepsilon}(t,x_i)+I_{\mathbf{N},\varepsilon}(t,x_i)+R_{\mathbf{N},\varepsilon}(t,x_i)\Big)=\eps^{-d} $, since this quantity does not depend upon $t$. It then follows that  $ 0\le Z_{\mathbf{N},\varepsilon}(t) \le \eps^{-d}$, for all $t\ge 0 $.  
Then by letting $ \dis \beta_{T}^{\mathbf{N},\varepsilon}= \sup\limits_{\substack{ 1\le j  \le \ell \\  0\le t \le T}} \beta_j\big(t,Z_{\mathbf{N},\varepsilon}(t)\big)$, we have that $ \dis \beta_{T}^{\mathbf{N},\varepsilon} \le  \overline{C}$, where $ \overline{C}= \max\left\{ \; \bar{\beta}, \bar{\alpha} , \dfrac{\mu_S}{\varepsilon^2} , \dfrac{\mu_I}{\varepsilon^2} , \dfrac{\mu_R}{\varepsilon^2}, \eps^{-d} \; \right\}. $
\begin{eqnarray}\label{fap}
Z_{\mathbf{N},\varepsilon}(t) 
&=& Z_{\mathbf{N},\varepsilon}(0) + \sum_{j=1}^{k_{\eps}}\frac{h_j}{\mathbf{N}} P_j\Bigg( \mathbf{N}\int_0^t \beta_j\big(r,Z_{\mathbf{N},\varepsilon}(r)\big)dr\Bigg).
\end{eqnarray}
Let $ \tau_1^j<\tau_2^j< \cdots $ be the jump times of the Poisson process $ P_j(t)$, $ 1 \le j\le k $ . As long as $ \mathbf{N} \beta_j\big( Z_{\mathbf{N},\varepsilon}(0)\big)\times~t~<~\tau_1^j $, for all $1\le j\le k$,  the process $Z_{\mathbf{N},\varepsilon}(t)$ remains constant. Let us set 
$$ T_1 = \inf\Big\{ t>0 : \mathbf{N}\beta_j\Big( Z_{\mathbf{N},\varepsilon}(0)\Big)\times t = \tau_1^j , \; \text{for some }  \; 1\le j\le k \Big\}. $$ 
The independence of the $P_j$'s ensures that there is almost surely a unique $j_{_0}$ such that \\ $ \mathbf{N}\beta_{j_{_0}}\big( Z_{\mathbf{N},\varepsilon}(0)\big)\times T_1=\tau_1^{j_{_0}}$. In this case $ \dis Z_{\mathbf{N},\varepsilon}(T_1)= Z_{\mathbf{N},\varepsilon}(0) + \dfrac{h_{j_{_0}}}{\mathbf{N}}$, and the process remains constant until the next jump of one of the $ P_j$'s. We wait for the next time for which one of the integrands
$$ \int_0^t \mathbf{N} \beta_j\Big( Z_{\mathbf{N},\varepsilon}(r)\Big)dr = \mathbf{N}\beta_j\Big( Z_{\mathbf{N},\varepsilon}(0)\Big)\times T_1+ \mathbf{N}\beta_j\Big( Z_{\mathbf{N},\varepsilon}(0)+ \dfrac{h_{j_{_0}}}{\mathbf{N}}\Big)\big(t-T_1\big) $$
will be equal  to the jump time of one of the $ P_j.$ We  continue this procedure . Since there are a finite number of $ P_j $ and  the rates $\beta_j$ are bounded, any time $t\in[0,T]$ is achieved after a finite number of operations as above. This construction shows  existence and uniqueness of the solution of (\ref{fap}).

\section{Law of large numbers ($\mathbf{N} \to \infty $, $\varepsilon$ being fixed)}
\rhead{Law of large numbers}
Recall that, from \ref{a1},   $  \dis  \int_D \big(\mathbf{s}(0,x)+\mathbf{i}(0,x)+\mathbf{r}(0,x)\big)dx=1$. Recall that in the stochastic model, we have a total of $\mathbf{N}\eps^{-d}$ individuals. At time $t=0$, each individual, independently of the others, is susceptible and located at  site $x_i$ with probability $\dis \int_{V_i}\!\!\mathbf{s}(0,x)dx$, infectious and located at  site $x_i$ with probability $ \dis  \int_{V_i}\!\! \mathbf{i}(0,x)dx$, removed and located at site $x_i$ with probability $\dis \int_{V_i} \!\!\mathbf{r}(0,x)dx$, $1\le i\le \eps^{-d}$.  It follows from the choice of the initial condition of the stochastic system, the law of large numbers and the definition (\ref{rn})  that for any $1\le i\le \eps^{-d}$, as $N\to \infty$, $ \dis S_{\mathbf{N},\eps}(0,x_i)\longrightarrow\eps^{-d}\int_{V_i}\!\!\mathbf{s}(0,x)dx$,  $ \dis I_{\mathbf{N},\eps}(0,x_i)\longrightarrow\eps^{-d}\int_{V_i}\!\mathbf{i}(0,x)dx$ and  $ \dis R_{\mathbf{N},\eps}(0,x_i)\longrightarrow\eps^{-d}\int_{V_i}\!\!\mathbf{r}(0,x)dx$,   a.s. .

\vspace{0.3cm}

In this section we fix the mesh  size $ \varepsilon$  of the grid and we let $\mathbf{N}$ go to infinity. We will show that the stochastic model converges to the corresponding deterministic model on the grid. First let us recall the law of large numbers for  Poisson processes.
\begin{lem}\label{llnp}
	Let $ \dis \left\{ \; P(t), \; t\ge 0 \; \right\}$ be a rate $\lambda $ Poisson process. Then 
	$$ \dfrac{P(t)}{t} \longrightarrow \lambda \; \text{a.s} \quad \text{as} \; t\to \infty .$$
\end{lem}
A proof of this well-known lemma can be found e.g. in Britton and Pardoux (2019).\fpr 
In the sequel, we shall assume that $Z_{\mathbf{N},\ep}(t)$ is defined on the probability space $\dis \Big(\Omega, \mathcal{F}, \mathcal{F}_t^{\mathbf{N},\ep}, \P\Big)$, where  $ \dis  \mathcal{F}_t^{\mathbf{N},\varepsilon} = \sigma\{ Z_{\mathbf{N},\varepsilon}(r,x_i) , \;  0\leq r \leq t \, ; \;  x_i \in D_{\varepsilon} \}$. If we consider the $k_\ep$-dimensional process $\dis\Big(M_{\mathbf{N},\varepsilon}^j\Big)_{1\le j \le k_\ep} $  whose  j-th  component is defined as
$\dis M_{\mathbf{N},\varepsilon}^j(t,x_i) = \dfrac{1}{\mathbf{N}}\mathrm{M}_j\Big( \mathbf{N}\int_0^t \beta_j\big(Z_{\mathbf{N},\varepsilon}(r,x_i)\big)dr\Big),$ for a site  $x_i\in  D_{\varepsilon}$, then we have the following Proposition.
\begin{prop}\label{normum}  
	For all $ 1\le j\le k_\ep$ for all $T>0$, as $ \mathbf{N}\rightarrow +\infty $ , 
	\[  \underset{0\leq t\leq T}{\sup}{\Big\rvert M_{\mathbf{N},\varepsilon}^j(t,x_i)\Big\rvert}\overset{a.s.}{\longrightarrow}0 .\]
\end{prop}
\begin{pr}  
	For all $T>0$ we have
	\begin{eqnarray}
	\underset{0\leq t\leq T}{\sup}{\Big\rvert M_{\mathbf{N},\varepsilon}^j(t,x_i)\Big\rvert}
	& = & \underset{0\leq t\leq T}{\sup}\Bigg\rvert \frac{1}{\mathbf{N}}M_j\left(\int_0^t \mathbf{N}\beta_j\big(Z_{\mathbf{N},\varepsilon}(r,x_i)\big)dr \right) \Bigg\rvert \nonumber\\
	& \leq & \underset{0\leq t\leq T\overline{C}}{\sup}\Big\rvert \frac{1}{\mathbf{N}}M_j(\mathbf{N}t) \Big\rvert \qquad \Big(\, \text{because} \; 0 \le \beta_j \le \overline{C} \,  \Big)  \nonumber \\
	& = & \underset{0\leq t\leq T\overline{C}}{\sup}\Big\rvert \frac{1}{\mathbf{N}}P_j(\mathbf{N}t)- t \Big\rvert. \nonumber
	\end{eqnarray}
	From \ref{llnp},  $$ \frac{P_j(\mathbf{N}\,t)}{\mathbf{N}}\longrightarrow t \; \; \text{a.s.}\, ,\; \text{as}\; \mathbf{N}\rightarrow \infty . $$
	We have pointwise convergence of a sequence of increasing functions towards a continous function, then from  the second  Dini Theorem this convergence is uniform on any compact time interval. This shows that 
	$$ \underset{0\leq t\leq T\bar{C}}{\sup}\Big\rvert \frac{1}{\mathbf{N}}P_j(\mathbf{N}\,t)- t \Big\rvert \longrightarrow 0 \; \; \text{a.s.}\, , \; \;  \text{as}\; \mathbf{N}\,\rightarrow \infty $$ and the Proposition is established.
	\fpr
\end{pr}
In what follows, $\Vert u \Vert $ denotes the norm  of an $\ell$–dimensional vector u defined as follow $\dis \Vert u \Vert := \sum_{j=1}^{\ell} \vert u_j \vert $. \\ Now we can prove the main result of this section. This law of large numbers is in fact a particular case of the general result in Britton $\&$ Pardoux (2019). But since the proof is rather short, we prefered to include it for the convenience of the reader.
\begin{thm} \label{llns} $\bf{(Law \ of \ Large \ Numbers)}$
	\rhead{Law \ of \ Large \ Numbers}
	\\Let  $ Z_{\mathbf{N},\varepsilon} $ denote the  solution  of the SDE (\ref{eqdifs})  and $ Z_{\varepsilon}$ the solution of the ODE \; $\dfrac{dZ_{\varepsilon}(t)}{dt}=~ b_{\varepsilon}(t,Z_{\varepsilon}(t))$. \\
	Let us fix an arbitrary $T > 0$ and assume that $\Big\Vert Z_{\mathbf{N},\varepsilon}(0)-Z_{\varepsilon}(0)\Big\Vert \longrightarrow 0 $,  as $ \mathbf{N}\rightarrow + \infty.$\\
	Then $ \dis   \underset{0\leq t\leq T}{\sup}\Big\Vert  Z_{\mathbf{N},\varepsilon}(t)-Z_{\varepsilon}(t) \Big\Vert \longrightarrow 0  \; \text{a.s.} \; , \; \; as \; \;   \mathbf{N}\rightarrow + \infty $ . 
\end{thm}

\begin{pr} Let us define
	$\dis M_{\mathbf{N},\varepsilon}(t) = \sum_{j=1}^{k_\ep} h_j M_{\mathbf{N},\varepsilon}^j(t), \; t\in [0,T] .$ We  first note that 
	
	$\dis \hspace{4cm}	\underset{0\leq t \leq T}{\sup}\Big\Vert M_{\mathbf{N},\varepsilon}(t) \Big\Vert
	\leq \sum_{j=1}^{k_\ep}\Vert h_j \Vert \underset{0\leq t \leq T} {\sup}\Big\rvert M_{\mathbf{N},\varepsilon}^j(t) \Big\rvert.$
	
	Hence from \ref{normum}, we deduce that $\underset{0\leq t \leq T} {\sup}\Big\Vert M_{\mathbf{N},\varepsilon}(t) \Big\Vert \overset{a.s}{\longrightarrow} 0,$ as $\mathbf{N}\rightarrow +\infty .$ \\
	Next for any $ r \in [0,T] $ we have
	\begin{eqnarray}
	\bigg\Vert b_{\varepsilon}\big(r, Z_{\mathbf{N},\varepsilon}(r)\big) - b_{\varepsilon}\big(r, Z_{\varepsilon}(r)\big)\bigg\Vert \nonumber\\
	&& \hspace{-4cm} = 2  \sum_{i=1}^{\ell}\beta(x_i)\bigg\rvert \dfrac{S_{\mathbf{N},\varepsilon}(r,x_i)I_{\mathbf{N},\varepsilon}(r,x_i)}{S_{\mathbf{N},\varepsilon}(r,x_i)+I_{\mathbf{N},\varepsilon}(r,x_i)+R_{\mathbf{N},\varepsilon}(r,x_i)}- \dfrac{S_{\varepsilon}(r,x_i)I_{\varepsilon}(r,x_i)}{S_{\varepsilon}(r,x_i)+I_{\varepsilon}(r,x_i)+R_{\varepsilon}(r,x_i)}\bigg\rvert \nonumber \\
	& & \hspace{-4cm} + 2 \sum_{i=1}^{\ell}\alpha(x_i) \Big\rvert I_{\mathbf{N},\varepsilon}(r,x_i)- I_{\varepsilon}(r,x_i)\Big\rvert +\mu_S\sum_{i=1}^{\ell} \Big\rvert \Delta_{\varepsilon} \Big(S_{\mathbf{N},\varepsilon}(r,x_i)-S_{\varepsilon}(r,x_i)\Big)\Big\vert  \nonumber \\
	& &\hspace{-4cm}+\mu_I\sum_{i=1}^{\ell} \Big\rvert \Delta_{\varepsilon} \Big(I_{\mathbf{N},\varepsilon}(r,x_i)-I_{\varepsilon}(r,x_i)\Big)\Big\vert +\mu_R\sum_{i=1}^{\ell} \Big\rvert \Delta_{\varepsilon} \Big(R_{\mathbf{N},\varepsilon}(r,x_i)- R_{\varepsilon}(r,x_i)\Big)\Big\rvert. \nonumber
	\end{eqnarray}
	Then, the fact that $ \beta$  and $ \alpha $ are bounded leads to 
	\begin{eqnarray}
	\hspace{-1.5cm}\bigg\Vert b_{\varepsilon}\big(r, Z_{\mathbf{N},\varepsilon}(r)\big) - b_{\varepsilon}\big(r, Z_{\varepsilon}(r)\big)\bigg\Vert \nonumber\\
	&&\hspace{-4cm}\leq 2\bar{\beta}\sum_{i=1}^{\ell}\Bigg\{2\Big\rvert S_{\mathbf{N},\varepsilon}(r,x_i)- S_{\varepsilon}(r,x_i)\Big\rvert +2\Big\rvert I_{\mathbf{N},\varepsilon}(r,x_i)- I_{\varepsilon}(r,x_i)\Big\rvert+ \Big\rvert R_{\mathbf{N},\varepsilon}(r,x_i)- R_{\varepsilon}(r,x_i)\Big\rvert \Bigg\} \nonumber\\
	&&\hspace{-4cm} + 2\bar{\alpha} \sum_{i=1}^{\ell}\Big\rvert I_{\mathbf{N},\varepsilon}(r,x_i)- I_{\varepsilon}(r,x_i)\Big\rvert +4\mu_S\,\varepsilon^{-2}\sum_{i=1}^{\ell} \Big\rvert S_{\mathbf{N},\varepsilon}(r,x_i)-S_{\varepsilon}(r,x_i)\Big\vert \nonumber\\
	&&\hspace{-4cm}+4\mu_I\,\varepsilon^{-2}\sum_{i=1}^{\ell} \Big\rvert I_{\mathbf{N},\varepsilon}(r,x_i)- I_{\varepsilon}(r,x_i)\Big\rvert +4\mu_R\,\varepsilon^{-2}\sum_{i=1}^{\ell} \Big\rvert R_{\mathbf{N},\varepsilon}(r,x_i)- R_{\varepsilon}(r,x_i)\Big\rvert \nonumber\\
	&&\hspace{-4cm}  \le C(\bar{\alpha},\bar{\beta},\bar{\mu},\varepsilon)\,\Big\Vert Z_{\mathbf{N},\varepsilon}(r) - Z_{\varepsilon}(r)\Big\Vert, \; \; \text{where}\;   \bar{\mu}=\max\{\mu_S,\mu_I,\mu_R\}. \nonumber
	\end{eqnarray}
	
	Hence  we have for all $t\in [0,T]$
	\begin{align*}
	\Big\Vert Z_{\mathbf{N},\varepsilon}(t)- Z_{\varepsilon}(t)\Big\Vert 
	& \leq \Big\Vert Z_{\mathbf{N},\varepsilon}(0)- Z_{\varepsilon}(0)\Big\Vert + \int_0^t \Big\Vert b_{\varepsilon}\big(r,Z_{\mathbf{N},\varepsilon}(r)\big)- b_{\varepsilon}\big(r,Z_{\varepsilon}(r)\big)\Big\Vert dr +  \Big\Vert M_{\mathbf{N},\varepsilon}(t) \Big\Vert\\
	& \leq \Bigg( \Big\Vert Z_{\mathbf{N},\varepsilon}(0)- Z_{\varepsilon}(0) \Big\Vert + \Big\Vert M_{\mathbf{N},\varepsilon}(t) \Big\Vert \Bigg) + C(\bar{\alpha},\bar{\beta},\bar{\mu},\varepsilon) \int_0^t \Big\Vert Z_{\mathbf{N},\varepsilon}(r)- Z_{\varepsilon}(r)\Big\Vert dr,
	\end{align*}
	and it follows from Gronwall’s Lemma that
	\[ \underset{0\leq t\leq T}{\sup}\Big\Vert Z_{\mathbf{N},\varepsilon}(t)- Z_{\varepsilon}(t)\Big\Vert \leq \left( \Big\Vert  Z_{\mathbf{N},\varepsilon}(0) - Z_{\varepsilon}(0) \Big\Vert + \underset{0\leq t\leq T}{\sup}\Big\Vert M_{\mathbf{N},\varepsilon}(t) \Big\Vert\right)\exp\Big( C(\bar{\alpha},\bar{\beta},\bar{\mu},\varepsilon) T\Big) .\]  This concludes  the proof of the theorem , since 
	$\dis \Big\Vert  Z_{\mathbf{N},\varepsilon}(0) - Z_{\varepsilon}(0) \Big\Vert + \underset{0\leq t\leq T}{\sup}\Big\Vert M_{\mathbf{N},\varepsilon}(t) \Big\Vert \longrightarrow 0,  \; \text{as} \;   \mathbf{N} \rightarrow +\infty. $
	\begin{flushright}
		$\blacksquare$
	\end{flushright}
\end{pr} 
We have just shown that the solution of the stochastic model $(\ref{eqdifs})$ converges a.s. locally
uniformly in $t$ to the solution of the deterministic model (\ref{eqdet}), as $ \mathbf{N} \to \infty $, $ \varepsilon $  being fixed.
If we then let $ \varepsilon \to 0 $,  we know from \ref{ct} that the discrete deterministic system converges  in $ L^{\infty}(D) $  towards the system of PDEs on D
\begin{equation*}
\hspace{-4.5cm}\left \{
\begin{aligned}
\dfrac{\partial\,\mathbf{s}}{\partial t}(t,x)= & \;-\dfrac{\beta(x)\, \mathbf{s}(t,x)\mathbf{i}(t,x) }{\mathbf{s}(t,x)+\mathbf{i}(t,x)+\mathbf{r}(t,x)}  + \mu_S \,\Delta \mathbf{s}(t,x) \\
\dfrac{\partial\,\mathbf{i}}{\partial t}(t,x)=&\;\dfrac{\beta(x)\, \mathbf{s}(t,x)\mathbf{i}(t,x) }{\mathbf{s}(t,x)+\mathbf{i}(t,x)+\mathbf{r}(t,x)} - \alpha(x)\, \mathbf{i}(t,x)+\mu_I\, \Delta \mathbf{i}(t,x)\\
\dfrac{\partial\,\mathbf{r}}{\partial t}(t,x)=&\;\alpha(x)\, \mathbf{i}(t,x)+\mu_S\, \Delta \mathbf{r}(t,x),  \quad (t,x) \in (0,T)\times  D. 
\end{aligned}
\right.
\end{equation*}

\section{Law of Large Numbers in the Supremum norm}
\rhead{Law of Large Numbers in Sup-norm}

In this section we let both the population size go to infinity and the mesh size $\varepsilon$ of the grid go to zero. Under the weak condition $\dfrac{\mathbf{N}}{\log(1/\varepsilon)}\longrightarrow \infty$, we obtain that the stochastic spatial model converges in probability to the corresponding deterministic one.\\
Let us define $\dis \;   \mathcal{S}_{\mathbf{N},\varepsilon}(t,x) = \sum_{i=1}^{\varepsilon^{-d}} S_{\mathbf{N},\varepsilon}(t,x_i)\mathbf{1}_{V_i}(x),\; \;  $   
$\dis \mathcal{I}_{\mathbf{N},\varepsilon}(t,x) = \sum_{i=1}^{\varepsilon^{-d}} I_{\mathbf{N},\varepsilon}(t,x_i)\mathbf{1}_{V_i}(x), \; $ and \\
$\dis \mathcal{R}_{\mathbf{N},\varepsilon}(t,x) = \sum_{i=1}^{\varepsilon^{-d}} R_{\mathbf{N},\varepsilon}(t,x_i)\mathbf{1}_{V_i}(x)$, $ (t,x)\in [0,T]\times D. $\; 
$ \Big(\mathcal{S}_{\mathbf{N},\varepsilon},  \mathcal{I}_{\mathbf{N},\varepsilon}, \mathcal{R}_{\mathbf{N},\varepsilon} \Big)$ is  solution of the SDEs
\begin{equation}
\left\{ 
\begin{aligned} 
\bigskip
\mathcal{S}_{\mathbf{N},\varepsilon}(t,x)
& = \mathcal{S}_{\mathbf{N},\varepsilon}(0,x)  + \mu_S \int_0^t\!\! \Delta_{\varepsilon}\mathcal{S}_{\mathbf{N},\varepsilon}(r,x)dr - \int_0^t \! \dfrac{\beta_\ep(x) \mathcal{S}_{\mathbf{N},\varepsilon}(r,x)\mathcal{I}_{\mathbf{N},\varepsilon}(r,x)}{\mathcal{S}_{\mathbf{N},\varepsilon}(r,x)+\mathcal{I}_{\mathbf{N},\varepsilon}(r,x)+\mathcal{R}_{\mathbf{N},\varepsilon}(r,x)}dr\\
&\hspace{1cm} + \mathcal{M}_{\mathbf{N},\varepsilon}^S(t,x)\\
\bigskip 
\mathcal{I}_{\mathbf{N},\varepsilon}(t,x) & = \mathcal{I}_{\mathbf{N},\varepsilon}(0,x) + \mu_I \int_0^t\!\! \Delta_{\varepsilon}\mathcal{I}_{\mathbf{N},\varepsilon}(r,x)dr + \int_0^t\! \dfrac{\beta_\ep(x) \mathcal{S}_{\mathbf{N},\varepsilon}(r,x)\mathcal{I}_{\mathbf{N},\varepsilon}(r,x)}{\mathcal{S}_{\mathbf{N},\varepsilon}(r,x)+\mathcal{I}_{\mathbf{N},\varepsilon}(r,x)+\mathcal{R}_{\mathbf{N},\varepsilon}(r,x)}dr \\
&\hspace{1cm}- \int_0^t\!\alpha_{\ep}(x)\mathcal{I}_{\mathbf{N},\varepsilon}(r,x)dr + \mathcal{M}_{\mathbf{N},\varepsilon}^I(t,x) \\
\mathcal{R}_{\mathbf{N},\varepsilon}(t,x) & = \mathcal{R}_{\mathbf{N},\varepsilon}(0,x) + \mu_R \int_0^t\!\! \Delta_{\varepsilon}\mathcal{R}_{\mathbf{N},\varepsilon}(r,x)dr + \int_0^t \! \alpha_{\ep}(x)\mathcal{I}_{\mathbf{N},\varepsilon}(r,x)dr + \mathcal{M}_{\mathbf{N},\varepsilon}^R(t,x)\\
(t,x) &\in [0,T]\times  D ,
\end{aligned}
\right. 
\end{equation} 

\vspace{-0.7cm}

\begin{align*}
\text{where}\quad \qquad  \mathcal{M}_{\mathbf{N},\varepsilon}^S(t,x)& =-\frac{1}{\mathbf{N}}\sum_{i=1}^{\varepsilon^{-d}} \mathrm{M}_{x_i}^{inf}\left(\mathbf{N}\int_0^t \dfrac{\beta(x_i) S_{\mathbf{N},\varepsilon}(r,x_i)I_{\mathbf{N},\varepsilon}(r,x_i)}{S_{\mathbf{N},\varepsilon}(r,x_i)+I_{\mathbf{N},\varepsilon}(r,x_i)+R_{\mathbf{N},\varepsilon}(r,x_i)}dr \right)\mathbf{1}_{V_i}(x)\\
& - \frac{1}{\mathbf{N}}\sum_{i=1}^{\varepsilon^{-d}} \sum_{y_i\sim x_i}\mathrm{M}_{S,x_i,y_i}^{mig}\left(\dfrac{\mu_S \mathbf{N}}{\varepsilon^2} \int_0^t S_{\mathbf{N},\varepsilon}(r,x_i)dr \right)\mathbf{1}_{V_i}(x)  \\
& + \frac{1}{\mathbf{N}}\sum_{i=1}^{\varepsilon^{-d}} \sum_{y_i\sim x_i}\mathrm{M}_{S,y_i,x_i}^{mig}\left( \dfrac{\mu_S \mathbf{N}}{\varepsilon^2} \int_0^tS_{\mathbf{N}, \varepsilon}(r,y_i)dr \right)\mathbf{1}_{V_i}(x) ,
\end{align*}
\begin{align*}
\mathcal{M}_{\mathbf{N},\varepsilon}^I(t,x)& =\frac{1}{\mathbf{N}}\sum_{i=1}^{\varepsilon^{-d}} \mathrm{M}_{x_i}^{inf}\left(\mathbf{N}\int_0^t\dfrac{\beta(x_i) S_{\mathbf{N},\varepsilon}(r,x_i)I_{\mathbf{N},\varepsilon}(r,x_i)}{S_{\mathbf{N},\varepsilon}(r,x_i)+I_{\mathbf{N},\varepsilon}(r,x_i)+R_{\mathbf{N},\varepsilon}(r,x_i)}dr\right)\mathbf{1}_{V_i}(x)\\
& - \frac{1}{\mathbf{N}}\sum_{i=1}^{\varepsilon^{-d}}\mathrm{M}_{x_i}^{rec}\left(\mathbf{N}\int_0^t\alpha(x_i)I_{\mathbf{N},\varepsilon}(r,x_i)dr \right)\mathbf{1}_{V_i}(x)\\
& - \frac{1}{\mathbf{N}}\sum_{i=1}^{\varepsilon^{-d}} \sum_{y_i\sim x_i}\mathrm{M}_{I,x_i,y_i}^{mig}\left(\dfrac{\mu_I \mathbf{N}}{\varepsilon^2} \int_0^t I_{\mathbf{N},\varepsilon}(r,x_i)dr \right)\mathbf{1}_{V_i}(x)  \\
& + \frac{1}{\mathbf{N}}\sum_{i=1}^{\varepsilon^{-d}} \sum_{y_i\sim x_i}\mathrm{M}_{I,y_i,x_i}^{mig}\left( \dfrac{\mu_I\mathbf{N}}{\varepsilon^2} \int_0^tS_{\mathbf{N}, \varepsilon}(r,y_i)dr \right)\mathbf{1}_{V_i}(x) , 
\end{align*}

\vspace{-0.8cm}

\begin{eqnarray}
\hspace{-10cm}\mathcal{M}_{\mathbf{N},\varepsilon}^R(t,x)& = &- \frac{1}{\mathbf{N}}\sum_{i=1}^{\varepsilon^{-d}}\mathrm{M}_{x_i}^{rec}\left(\mathbf{N}\int_0^t\alpha(x_i)I_{\mathbf{N},\varepsilon}(r,x_i)dr \right)\mathbf{1}_{V_i}(x)\nonumber\\
& - &\frac{1}{\mathbf{N}}\sum_{i=1}^{\varepsilon^{-d}} \sum_{y_i\sim x_i}\mathrm{M}_{R,x_i,y_i}^{mig}\left(\dfrac{\mu_R \mathbf{N}}{\varepsilon^2} \int_0^t R_{\mathbf{N},\varepsilon}(r,x_i)dr \right)\mathbf{1}_{V_i}(x) \nonumber \\
& +& \frac{1}{\mathbf{N}}\sum_{i=1}^{\varepsilon^{-d}} \sum_{y_i\sim x_i}\mathrm{M}_{R,y_i,x_i}^{mig}\left(\dfrac{\mu_R\mathbf{N}}{\varepsilon^2} \int_0^t R_{\mathbf{N}, \varepsilon}(r,y_i)dr \right)\mathbf{1}_{V_i}(x).\nonumber
\end{eqnarray}

Here we set 
\; $ \dis  X_{\mathbf{N}, \varepsilon} = \left( \begin{array}{cl}
\mathcal{S}_{\mathbf{N},\varepsilon} \\
\mathcal{I}_{\mathbf{N},\varepsilon}\\
\mathcal{R}_{\mathbf{N},\varepsilon}
\end{array}\right) $ \; and \;   $ \dis  \mathcal{M}_{\mathbf{N}, \varepsilon} = \left( \begin{array}{cl}
\mathcal{M}_{\mathbf{N},\varepsilon}^S \\
\mathcal{M}_{\mathbf{N},\varepsilon}^I \\
\mathcal{M}_{\mathbf{N},\varepsilon}^R
\end{array}\right).$ 
Recall that $ \mathcal{X}_{\varepsilon}= \left( \begin{array}{cl}
\mathcal{S}_{\varepsilon} \\
\mathcal{I}_{\varepsilon}\\
\mathcal{R}_{\varepsilon}
\end{array}\right) $ and  $ X= \left( \begin{array}{cl}
\mathbf{s}\\
\mathbf{i}\\
\mathbf{r}
\end{array}\right) $ .
\\ 

The main goal of this section is to  prove the following result.

\begin{thm}[\textbf{Law of Large Numbers in Sup-norm}]\label{llnsup}~ \\
	Let us assume that $(\ep,\mathbf{N})\to (0,\infty)$, in such way that
	\begin{enumerate}
		\item[(i)] $\dfrac{\mathbf{N}}{\log(1/\varepsilon)}\longrightarrow \infty$ as $\mathbf{N} \to \infty $ and $ \varepsilon\to 0 $;
		\item[(ii)] $\Big\Vert X_{\mathbf{N},\varepsilon}(0)- X(0) \Big\Vert_{\infty} \longrightarrow  0 $ in probability.
	\end{enumerate}	
	Then for all $ T > 0 $,  $ \un{t\in [0,T]}{\sup}\BV X_{\mathbf{N},\varepsilon}(t)- X(t) \BV_{\infty} \longrightarrow 0 $ in probability .
\end{thm}
~\\
We prove the Theorem in the case $d=2$, but the result holds true in dimensions d= 1, 3 as well, as we will explain below. \\
Since $\un{t\in [0,T]}{\sup}\BV \mathcal{X}_{\varepsilon}(t)- X(t) \BV_{\infty} \longrightarrow 0 $ by \ref{ct}, clearly our Theorem will follow from 

\begin{prop}\label{pe} For all $T>0$,
	$ \dis \un{t\in [0,T]}{\sup}\BV X_{\mathbf{N},\varepsilon}(t)- \mathcal{X}_{\varepsilon}(t) \BV_{\infty} \longrightarrow 0 \; \text{in probability} $,
	as $\mathbf{N}\to \infty $ and $\varepsilon\to 0 $, in such a way that $ \dfrac{\mathbf{N}}{\log(1/\varepsilon) } \longrightarrow  \infty. $ 
\end{prop}


\begin{pr} For all $t\in [0,T]$, We have
	
	$\dis \hspace{3cm} X_{\mathbf{N},\varepsilon}(t)= X_{\mathbf{N},\varepsilon}(0) +  \int_{0}^t \widetilde{\Delta}_{\ep}X_{\mathbf{N}, \varepsilon}(r)dr + \int_{0}^tG\Big(X_{\mathbf{N},\varepsilon}(r)\Big)dr + \mathcal{M}_{\mathbf{N},\varepsilon}(t), $ 
	
	$\dis \hspace{3cm} \mcx(t)= \mathcal{X}_{\varepsilon}(0) + \int_{0}^t\widetilde{\Delta}_{\ep}\mcx(r)dr + \int_{0}^tG\Big(\mcx(r)\Big)dr, $
	
	\vspace{-0.6cm}
	
	\begin{eqnarray}
	X_{\mathbf{N},\varepsilon}(t)-\mcx(t)
	&=& \widetilde{\mathsf{T}}_{\ep}(t)\big[X_{\mathbf{N},\varepsilon}(0)-\mcx(0) \big] +\int_0^t \widetilde{\mathsf{T}}_{\ep}(t-r)\Big[ G\big(X_{\mathbf{N}, \varepsilon}(r)\big)-G\big(\mcx(r)\big)\Big] dr + Y_{\mathbf{N},\varepsilon}(t), \nonumber
	\end{eqnarray}
	where again $ \dis
	\widetilde{\mathsf{T}}_{\ep}(t)X_{\mathbf{N},\varepsilon}= \left( 
	\begin{array}{cl}
	\mteps(t)\mathcal{S}_{\mathbf{N},\varepsilon}\\
	\mtepi(t)\mathcal{I}_{\mathbf{N},\varepsilon} \\
	\mtepr(t)\mathcal{R}_{\mathbf{N},\varepsilon} 
	\end{array} \right)  $ and similarly for $ \widetilde{\mathsf{T}}_{\ep}(t) \mathcal{X}_{\varepsilon},\cdots $; \; $ \dis
	Y_{\mathbf{N},\varepsilon}(t)= \left( 
	\begin{array}{cl}
	Y_{\mathbf{N},\varepsilon}^S(t) \\
	Y_{\mathbf{N},\varepsilon}^I(t)\\
	Y_{\mathbf{N},\varepsilon}^R(t) 
	\end{array} \right) 
	$  and \\ $ \dis Y_{\mathbf{N},\varepsilon}^S(t)= \int_0^t \mteps(t-r)d\mathcal{M}_{\mathbf{N},\varepsilon}^S(r), \; \; Y_{\mathbf{N},\varepsilon}^I(t)= \int_0^t \mtepi(t-r)d\mathcal{M}_{\mathbf{N},\varepsilon}^I(r)$, $\dis Y_{\mathbf{N},\varepsilon}^R(t)= \int_0^t \mtepr(t-r)d\mathcal{M}_{\mathbf{N},\varepsilon}^R(r) . $ 
	As in the proof of \ref{llns}, one can show that there is a constant $C(\bar{\beta},\bar{\alpha})$ such that for all $ r \in [0,T] $,  we have 
	\begin{eqnarray}\label{mg}
	\BV G\Big(X_{\mathbf{N},\varepsilon}(r)\Big)-G\Big(\mcx(r)\Big) \BV_{\infty} 
	&\le & C(\bar{\alpha}, \bar{\beta}) \BV X_{\mathbf{N}, \varepsilon}(r) - \mcx(r) \BV_{\infty},
	\end{eqnarray}
	since G is globally Lipschitz.
	Using (\ref{mg}) and the fact that $ \widetilde{\mathsf{T}}_{\ep} $ is a contraction semigroup in $ \big(L^{\infty}(D)\big)^3$, we have 
	\begin{eqnarray}
	\BV X_{\mathbf{N},\varepsilon}(t) -\mcx(t)\BV_{\infty}
	&\le & \BV X_{\mathbf{N},\varepsilon}(0) - \mcx(0)\BV_{\infty} + C(\bar{\alpha}, \bar{\beta}) \int_0^t \BV X_{\mathbf{N},\varepsilon}(r) -\mcx(r)\BV_{\infty} dr + \BV Y_{\mathbf{N},\varepsilon}(t)\BV_{\infty}. \nonumber 
	\end{eqnarray} 
	It then follows from Gronwall's Lemma that
	\begin{eqnarray}\label{gs}
	\hspace{-1.5cm}\un{t\in [0,T]}{\sup}\BV X_{\mathbf{N},\varepsilon}(t) -\mcx(t)\BV_{\infty} 
	&\le &  \Bigg( \BV X_{\mathbf{N},\varepsilon}(0) -\mcx(0)\BV_{\infty} + \un{t\in [0,T]}{\sup}\BV Y_{\mathbf{N},\varepsilon}(t)\BV_{\infty} \Bigg) e^{C(\bar{\alpha}, \bar{\beta})T} .
	\end{eqnarray}
	Since
	$ \BV X_{\mathbf{N},\varepsilon}(0) -\mcx(0)\BV_{\infty} \longrightarrow 0 \; \; \text{in probability}, $ the Proposition follows from (\ref{gs}) and  \ref{ny} below.
	\fpr
\end{pr}
\begin{prop}\label{ny} For all $T>0$
	\begin{eqnarray}
	\un{t\in [0,T]}{\sup}\BV Y_{\mathbf{N},\varepsilon}(t)\BV_{\infty} \longrightarrow 0 \; \text{in probability} , \; \text{as $\mathbf{N}\to \infty $ and $\varepsilon\to 0 $, provided $ \dfrac{\mathbf{N}}{\log(1/\varepsilon) } \longrightarrow  \infty. $}
	\end{eqnarray}
\end{prop}

\vspace{0.4cm}

Before proving this Proposition, we first establish some technical Lemmas.

\begin{lem}\label{dm}
	Let $ \dis f= \varepsilon^{-2}\mathbf{1}_{V_i}$. Then, for any $J\in \{S, I, R\}$
	$$ \big\langle \;  \big(\nabla_{\varepsilon}^{1,+}\mtepj(t)f\big)^2 + \big(\nabla_{\varepsilon}^{1,-}\mtepj(t)f\big)^2 +\big(\nabla_{\varepsilon}^{2,+}\mtepj(t)f\big)^2 + \big(\nabla_{\varepsilon}^{2,-}\mtepj(t)f\big)^2 + \big(\mtepj(t)f\big)^2 , \; 1 \;  \big\rangle  \le h_{\varepsilon}(t) $$
	where
	\begin{eqnarray}\label{h}
	\int_0^t h_{\varepsilon}(r)dr\le C\,\varepsilon^{-2} + t .
	\end{eqnarray} 
\end{lem}

\begin{pr} For $ \dis f= \varepsilon^{-2}\mathbf{1}_{V_i}$ and $J\in \{S, I, R\}$, we have 
	\begin{eqnarray}
	\big\langle \;  \Big(\nabla_{\varepsilon}^{1,+}\mtepj(t)f\Big)^2 + \Big(\nabla_{\varepsilon}^{2,+}\mtepj(t)f\Big)^2 \; , \; 1 \; \big\rangle
	\!\!\!\!&=& \!\!\!\! \big\langle \;  \nabla_{\varepsilon}^{1,+}\mtepj(t)f \; , \; \nabla_{\varepsilon}^{1,+}\mtepj(t)f \; \big\rangle  + \big\langle \;  \nabla_{\varepsilon}^{2,+}\mtepj(t)f , \nabla_{\varepsilon}^{2,+}\mtepj(t)f \; \big\rangle  \nonumber \\
	\!\!\!\!&=&\!\!\!\! - \big\langle \;  \nabla_{\varepsilon}^{1,-} \nabla_{\varepsilon}^{1,+}\mtepj(t)f , \mtepj(t)f \; \big\rangle - \big\langle \;  \nabla_{\varepsilon}^{2,-}\nabla_{\varepsilon}^{2,+}\mtepj(t)f , \mtepj(t)f \; \big\rangle  \nonumber \\
	\!\!\!\!&=& \!\!\!\!- \big\langle \; \Delta_{\varepsilon}\mtepj(t)f , \mtepj(t)(t)f \; \big\rangle.   \nonumber 
	\end{eqnarray}
	Using  the facts that $ \Delta_{\varepsilon}\mtepj(t)f=\mtepj(t)\Delta_{\varepsilon}f$, $\mtepj(t)$ is self-adjoint and  (\ref{rept}), we obtain  
	
	\begin{eqnarray}
	\big\langle \;  \Big(\nabla_{\varepsilon}^{1,+}\mtepj(t)f\Big)^2 +    \Big(\nabla_{\varepsilon}^{2,+}\mtepj(t)f\Big)^2 \; , \; 1 \; \big\rangle
	&=& - \big\langle \; \mtepj(t)\Delta_{\varepsilon}f , \mtepj(t)f \; \big\rangle  \nonumber \\
	&=&- \big\langle \; \Delta_{\varepsilon}f , \mtepj(2t)f \; \big\rangle  \nonumber \\
	&=&  \sum_{m} \big\langle \;  f,\mathbf{f}_{m}^{\varepsilon}\; \big\rangle^2 e^{-2\lambda_{m,J}^{\varepsilon}t}\lambda_{m,J}^{\varepsilon}\nonumber \\
	&\le & 4\sum_{m}e^{-2 \lambda_{m,J}^{\varepsilon}t}\lambda_{m,J}^{\varepsilon}. \nonumber
	\end{eqnarray}
	
	\vspace{-0.2cm}
	
	Similarly $\dis \big\langle \;  \big(\nabla_{\varepsilon}^{1,-}\mtepj(t)f\big)^2 + \big(\nabla_{\varepsilon}^{2,-}\mtepj(t)f\big)^2 , 1 \; \big\rangle  \le 4\sum_{m}e^{-2 \lambda_{m,J}^{\varepsilon}t}\lambda_{m,J}^{\varepsilon}$.
	Moreover, we have
	\begin{eqnarray}
	\big\langle \;  \big(\mtepj(t)f\big)^2 , 1 \; \big\rangle 
	&=& \big\langle \;  \mtepj(2t)f , f \; \big\rangle \nonumber \\
	&=& 1+ \sum_{m\ne (0,0)}e^{-2\lambda_{m,J}^{\varepsilon}t} \big\langle \;  f,\mathbf{f}_{m}^{\varepsilon}\; \big\rangle^2 \nonumber \\
	&\le& 1+ 4 \sum_{m\ne (0,0)}e^{-2\lambda_{m,J}^{\varepsilon}t}.\nonumber
	\end{eqnarray}
	
	\vspace{-0.2cm}
	
	So, the result holds with $\dis h_{\varepsilon}(t) = 1+ 8\sum_{m\ne (0,0)}e^{-2 \lambda_{m,J}^{\varepsilon}t}\big(\lambda_{m,J}^{\varepsilon}+1\big).$ 
	
	\vspace{-0.7cm}
	
	\begin{eqnarray}
	\int_0^th_{\varepsilon}(r)dr
	&\le &  t + 8 \sum_{m\ne(0,0)} \int_0^{+\infty} e^{-2\lambda_{m,J}^{\varepsilon}r}\big(1+\lambda_{m,J}^{\varepsilon}\big)dr \nonumber \\
	&=& t + 4 \sum_{m\ne(0,0)} \bigg( 1 + \dfrac{1}{\lambda_{m,J}^{\varepsilon}}\bigg)\nonumber \\
	&\le & t + 4 \sum\limits_{\substack{1\le m_1 \le \varepsilon^{-1}\\  1\le m_2 \le \varepsilon^{-1}\\ }} \Bigg(1 +\dfrac{1}{ c \, (m_1^2+m_2^2)}\Bigg)\nonumber \\
	& \le & t+ 4 \varepsilon^{-2} + \sum\limits_{\substack{ 1\le m_1 \le \varepsilon^{-1}\\  1\le m_2 \le \varepsilon^{-1}}} \dfrac{2}{c} \nonumber \\
	&=& t + \varepsilon^{-2} \Big(4 + \dfrac{2}{ c} \Big). \nonumber 
	\end{eqnarray}
\end{pr}
\vspace{-1cm}
\begin{flushright}
	$\blacksquare$
\end{flushright}
For any càdlàg process $ Z $, let  $ \delta Z(t)=  Z(t)- Z(t^-) $  denote its jump at time $t$.\\

We shall need below the
\begin{lem}\label{em}
	Let $M_{t}$ be a bounded martingale of finite variation defined on $[t_0 ,  t_1]$ with $M_{t_0}=0$ and satisfying
	\begin{enumerate}
		\item[(i)] $M$ is right-continuous with left limits,
		\item[(ii)] $  \vert\delta M_{t} \vert \le c $ \;  for  $t_0 \le t \le t_1$, where $c$ is a positive constant,
		\item[(iii)] $ \dis \sum_{t_0 \le s \le t } \big(\delta M_{s}\big)^2-\int_{t_0}^t h(s)ds$ is a  supermartingale, where $ h $ is a positive deterministic function.
	\end{enumerate} 
	Then  \;  $ \dis \E\Big( \exp\big(M_{t_1}\big)\Big)\le \exp\Big(\dfrac{e^c}{2}\int_{t_0}^{t_1}h(s)ds\Big). $ 
\end{lem}

\begin{pr}
	
	Let $ f(x)=e^x$. We have $ 0 \le f^{''}(x+y)=f(x)f(y)\le e^c f(x) $ \;  for all $y\le c $. \\
	For $ t_0\le t\le t_1 $ 
	
	\vspace{-0.9cm}
	
	\begin{align*}
	f(M_t) 
	&=1 + \int_{t_0}^{t}f^{\prime}(M_{s^-})dM_s+\sum_{t_0\le s\le t}\Big(f(M_s)- f(M_{s^-})-f^{\prime}(M_{s^-})\delta M_s  \Big)\\
	& \le 1 + \int_{t_0}^{t}f^{\prime}(M_{s^-})dM_s+\dfrac{e^c}{2}\sum_{t_0\le s\le t}f(M_{s^-})(\delta M_s)^2 \\
	&= 1 + \int_{t_0}^{t}f^{\prime}(M_{s^-})dM_s+\dfrac{e^c}{2}\Big(\sum_{t_0\le s\le t}f(M_{s^-})(\delta M_s)^2 - \int_{t_0}^{t}f(M_{s})h(s)ds \Big) + \dfrac{e^c}{2}\int_{t_0}^{t}f(M_{s})h(s)ds 
	\end{align*}
	where we used Taylor's formula and (ii) for the inequality.\\
	From (iii) and the fact that $ \dis \int_{t_0}^{t}f^{\prime}(M_{s^-})dM_s $ has mean zero, we deduce 
	
	\vspace{-0.7cm}
	
	\begin{align*}
	\E\Big(f(M_t)\Big)
	& \le 1 + \dfrac{e^c}{2}\int_{t_0}^{t}\E\Big( f(M_{s})\Big)h(s)ds.
	\end{align*}
	The result now follows from  Gronwall's inequality.
	\fpr
\end{pr}
\vspace{-0.6cm}

\begin{lem}\label{mart1}
	
	For any site  $ x_i\in D_{\varepsilon} $, the following are $ \mathcal{F}_t^{\mathbf{N},\varepsilon} $ mean zero martingales:
	\begin{eqnarray}\label{s}
	\sum_{0\le r\le t }\bigg(\delta S_{\mathbf{N},\varepsilon}(r,x_i) \bigg)^2- \dfrac{1}{\mathbf{N}}\int_{0}^{t}\dfrac{\beta(x_i)S_{\mathbf{N},\varepsilon}(r,x_i)I_{\mathbf{N},\varepsilon}(r,x_i)}{S_{\mathbf{N},\varepsilon}(r,x_i)+I_{\mathbf{N},\varepsilon}(r,x_i)+R_{\mathbf{N},\varepsilon}(r,x_i)}dr 
	&& \\
	&& \hspace{-11cm}-\dfrac{4\mu_S}{\mathbf{N}\varepsilon^2} \int_{0}^{t}S_{\mathbf{N},\varepsilon}(r,x_i)dr -\dfrac{\mu_S}{\mathbf{N}\varepsilon^2}\int_{0}^{t} \bigg(\sum_{j=1}^{2} S_{\mathbf{N},\varepsilon}(r,x_i +\varepsilon e_j) + \sum_{j=1}^{2} S_{\mathbf{N},\varepsilon}(r, x_i -\varepsilon e_j) \bigg)dr  \nonumber 
	\end{eqnarray}
	
	\vspace{-0.7cm}
	
	\begin{eqnarray}
	\hspace{1cm}\sum_{0\le r\le t}\bigg(\delta  I_{\mathbf{\mathbf{N}},\varepsilon}(r,x_i) \bigg)^2- \dfrac{1}{\mathbf{\mathbf{N}}}\int_{0}^{t}\dfrac{\beta(x_i)S_{\mathbf{N},\varepsilon}(r,x_i)I_{\mathbf{N},\varepsilon}(r,x_i)}{S_{\mathbf{N},\varepsilon}(r,x_i)+I_{\mathbf{N},\varepsilon}(r,x_i)+R_{\mathbf{N},\varepsilon}(r,x_i)}dr  - \dfrac{1}{\mathbf{N}} \int_{0}^{t}\alpha(x_i)I_{\mathbf{N},\varepsilon}(r,x_i)dr
	&&  \\
	&& \hspace{-14cm}-\dfrac{4\mu_I}{\mathbf{N}\varepsilon^2} \int_{0}^{t}I_{\mathbf{N},\varepsilon}(r,x_i)dr -\dfrac{\mu_I}{\mathbf{N}\varepsilon^2}\int_{0}^{t} \bigg(\sum_{j=1}^{2} I_{\mathbf{N},\varepsilon}(r,x_i +\varepsilon e_j) + \sum_{j=1}^{2} I_{\mathbf{N},\varepsilon}(r,x_i-\varepsilon e_j) \bigg)dr \nonumber 
	\end{eqnarray}
	\begin{eqnarray}
	\sum_{0\le r\le t }\bigg(\delta R_{\mathbf{\mathbf{N}},\varepsilon}(r,x_i) \bigg)^2- \dfrac{1}{\mathbf{N}} \int_{0}^{t}\alpha(x_i)I_{\mathbf{N},\varepsilon}(r,x_i)dr -\dfrac{4\mu_R}{\mathbf{N}\varepsilon^2} \int_{0}^{t}R_{\mathbf{N},\varepsilon}(r,x_i)dr
	&& \\
	&&\hspace{-9cm} -\dfrac{\mu_R}{\mathbf{N}\varepsilon^2}\int_{0}^{t} \bigg(\sum_{j=1}^{2} R_{\mathbf{N},\varepsilon}(r,x_i +\varepsilon e_j) + \sum_{j=1}^{2} R_{\mathbf{N},\varepsilon}(r,x_i-\varepsilon e_j) \bigg)dr \nonumber
	\end{eqnarray}
	
	\vspace{-0.9cm}
	
	\begin{eqnarray}\label{sd}
	\sum_{0\le r\le t}\bigg(\delta  S_{\mathbf{N},\varepsilon}(r,x_i) \bigg)\bigg(\delta  S_{\mathbf{N},\varepsilon}(r ,x_i \pm \varepsilon e_j) \bigg)+\dfrac{\mu_S}{\mathbf{N}\varepsilon^2} \int_{0}^{t}\Big( S_{\mathbf{N},\varepsilon}(r,x_i)+S_{\mathbf{N},\varepsilon}(r,x_i  \pm \varepsilon e_j)\Big)dr 
	\end{eqnarray}
	
	\vspace{-0.9cm}
	
	\begin{eqnarray}
	\sum_{0\le r\le t }\bigg(\delta  I_{\mathbf{N},\varepsilon}(r,x_i) \bigg)\bigg(\delta  I_{\mathbf{N},\varepsilon}(r ,x_i \pm \varepsilon e_j) \bigg)+\dfrac{\mu_I}{\mathbf{N}\varepsilon^2} \int_{0}^{t}\Big( I_{\mathbf{N},\varepsilon}(r,x_i)dr+I_{\mathbf{N},\varepsilon}(r,x_i \pm \varepsilon e_j)\Big)dr
	\end{eqnarray}
	
	\vspace{-0.8cm}
	
	\begin{eqnarray}
	\sum_{0\le r\le t }\bigg(\delta  R_{\mathbf{N},\varepsilon}(r,x_i) \bigg)\bigg(\delta  R_{\mathbf{N},\varepsilon}(r ,x_i \pm \varepsilon e_j) \bigg)+\dfrac{\mu_R}{\mathbf{N}\varepsilon^2} \int_{0}^{t}\Big( R_{\mathbf{N},\varepsilon}(r,x_i)dr+R_{\mathbf{N},\varepsilon}(r,x_i \pm \varepsilon e_j)\Big)dr 
	\end{eqnarray}
	$j=1 , 2  $.
\end{lem}
\begin{pr}
	The proof of this Lemma is based on the computation of the jumps. For  (\ref{s}), 
	we have
	\begin{eqnarray}
	\sum_{0\le r\le t}\bigg(\delta  S_{\mathbf{N},\varepsilon}(r,x_i) \bigg)^2
	&=& \frac{1}{\mathbf{N}^2}\mathrm{P}_{x_i}^{inf}\left(\mathbf{N}\int_0^{t} \dfrac{\beta(x_i)S_{\mathbf{N},\varepsilon}(r,x_i)I_{\mathbf{N},\varepsilon}(r,x_i)}{S_{\mathbf{N},\varepsilon}(r,x_i)+I_{\mathbf{\mathbf{N}},\varepsilon}(r,x_i)+R_{\mathbf{\mathbf{N}},\varepsilon}(r,x_i)}dr \right) \nonumber \\
	&+& \sum_{y_i\sim x_i}\frac{1}{\mathbf{N}^2}\mathrm{P}_{S,x_i,y_i}^{mig}\left(\frac{\mu_S \mathbf{N}}{\varepsilon^2} \int_0^{t}S_{\mathbf{N},\varepsilon}(r,x_i)dr \right)\nonumber \\
	&+& \sum_{y_i\sim x_i}\frac{1}{\mathbf{N}^2}\mathrm{P}_{S,y_i,x_i}^{mig}\left(\frac{\mu_S \mathbf{N}}{\varepsilon^2} \int_0^{t}S_{\mathbf{N},\varepsilon}(r,y_i)dr \right).\nonumber 
	\end{eqnarray}
	By writing each Poisson process as  $ M(t)+t $ , we then have
	\begin{eqnarray}
	\sum_{0\le r\le t }\big(\delta  S_{\mathbf{N},\varepsilon}(r,x_i) \big)^2- \dfrac{1}{\mathbf{N}}\int_{0}^{t}\dfrac{\beta(x_i)S_{\mathbf{N},\varepsilon}(r,x_i)I_{\mathbf{N},\varepsilon}(r,x_i)}{S_{\mathbf{N},\varepsilon}(r,x_i)+I_{\mathbf{\mathbf{N}},\varepsilon}(r,x_i)+R_{\mathbf{\mathbf{N}},\varepsilon}(r,x_i)}dr \nonumber\\
	&&\hspace{-12cm}-\dfrac{4\mu_S}{\mathbf{N}\varepsilon^2}\int_{0}^{t}S_{\mathbf{N},\varepsilon}(r,x_i)dr-\dfrac{\mu_S}{\mathbf{N}\varepsilon^2}\int_{0}^{t} \bigg(\sum_{j=1}^{2} S_{\mathbf{N},\varepsilon}(r,x_i +\varepsilon e_j) + \sum_{j=1}^{2} S_{\mathbf{N},\varepsilon}(r, x_i -\varepsilon e_j) \bigg)dr \nonumber\\
	&&\hspace{-12cm}=\frac{1}{\mathbf{N}^2}\mathrm{M}_{x_i}^{inf}\left(\mathbf{N}\int_0^{t}\dfrac{\beta(x_i)S_{\mathbf{N},\varepsilon}(r,x_i)I_{\mathbf{N},\varepsilon}(r,x_i)}{S_{\mathbf{N},\varepsilon}(r,x_i)+I_{\mathbf{\mathbf{N}},\varepsilon}(r,x_i)+R_{\mathbf{\mathbf{N}},\varepsilon}(r,x_i)}dr \right)+ \sum_{y_i\sim x_i}\frac{1}{\mathbf{N}^2}\mathrm{M}_{S,x_i,y_i}^{mig}\left(\frac{\mu_S \mathbf{N}}{\varepsilon^2} \int_0^{t }S_{\mathbf{N},\varepsilon}(r,x_i)dr \right)\nonumber \\
	&&\hspace{-12cm} +\sum_{y_i\sim x_i}\frac{1}{\mathbf{N}^2}\mathrm{M}_{S,y_i,x_i}^{mig}\left(\frac{\mu_S \mathbf{N}}{\varepsilon^2} \int_0^{t}S_{\mathbf{N},\varepsilon}(r,y_i)dr \right), \nonumber
	\end{eqnarray}
	which is a martingale. The other statements are proved similarly.
	\fpr
\end{pr}
The following result is a consequence of the previous Lemma.

\begin{lem}\label{mmj}
	Let $ \varphi\in H^{\varepsilon}$. The following are mean zero  martingales
	\begin{eqnarray}\label{ms}
	\sum_{0\le r\le t}\bigg(\delta \langle \; \mathcal{M}_{\mathbf{N},\varepsilon}^S(r) \, , \, \varphi \, \rangle \bigg)^2 -\dfrac{\varepsilon^2}{\mathbf{N}} \int_0^{t}\langle \; \dfrac{\beta_{\ep}(.)\mathcal{S}_{\mathbf{N},\varepsilon}(r)\mathcal{I}_{\mathbf{N},\varepsilon}(r)}{\mathcal{S}_{\mathbf{N},\varepsilon}(r)+\mathcal{I}_{\mathbf{\mathbf{N}},\varepsilon}(r)+\mathcal{R}_{\mathbf{\mathbf{N}},\varepsilon}(r)}\, , \,  \varphi^2  \; \rangle \, dr\\
	&& \hspace{-12cm}-\dfrac{\mu_S\varepsilon^2}{\mathbf{N}} \int_0^{t}\langle\; \mathcal{S}_{\mathbf{N},\varepsilon}(r)\, ,\,  \big(\nabla_{\varepsilon}^{1,+}\varphi \big)^2+\big( \nabla_{\varepsilon}^{1,-}\varphi \big)^2+\big( \nabla_{\varepsilon}^{2,+}\varphi \big)^2+\big( \nabla_{\varepsilon}^{2,-}\varphi \big)^2 \; \rangle \, dr  \nonumber 
	\end{eqnarray}
	\begin{eqnarray}\label{mi}
	\hspace{1cm}\sum_{0\le r\le t}\bigg(\delta \langle \; \mathcal{M}_{\mathbf{N},\varepsilon}^I(r) \, , \, \varphi \, \rangle \bigg)^2 -\dfrac{\varepsilon^2}{\mathbf{N}} \int_0^{t}\langle \; \dfrac{\beta_{\ep}(.)\mathcal{S}_{\mathbf{N},\varepsilon}(r)\mathcal{I}_{\mathbf{N},\varepsilon}(r)}{\mathcal{S}_{\mathbf{N},\varepsilon}(r)+\mathcal{I}_{\mathbf{\mathbf{N}},\varepsilon}(r)+\mathcal{R}_{\mathbf{\mathbf{N}},\varepsilon}(r)}\, , \,  \varphi^2  \; \rangle dr \\
	&&\hspace{-14cm}- \dfrac{ \varepsilon^2}{\mathbf{N}}\int_0^{t} \langle\; \alpha_{\ep}(.)\, \mathcal{I}_{\mathbf{N}, \varepsilon}(r) , \varphi^2 \; \rangle dr
	-\dfrac{\mu_I \,\varepsilon^2}{\mathbf{N}} \int_0^{t}\langle\; \mathcal{I}_{\mathbf{N},\varepsilon}(r)\, ,\,  \big( \nabla_{\varepsilon}^{1,+}\varphi \big)^2+\big( \nabla_{\varepsilon}^{1,-}\varphi \big)^2+\big( \nabla_{\varepsilon}^{2,+}\varphi \big)^2+\big( \nabla_{\varepsilon}^{2,-}\varphi \big)^2 \; \rangle \, dr  \nonumber  
	\end{eqnarray}
	\begin{eqnarray}\label{mr}
	\sum_{0\le r\le t }\bigg(\delta \langle \; \mathcal{M}_{\mathbf{N},\varepsilon}^R(r) \, , \, \varphi \, \rangle \bigg)^2  - \dfrac{ \varepsilon^2}{\mathbf{N}}\int_0^{t} \langle\; \alpha_{\ep}(.)\, \mathcal{I}_{\mathbf{N}, \varepsilon}(r) , \varphi^2 \; \rangle dr \\
	&& \hspace{-9cm}-\dfrac{\mu_R \,\varepsilon^2}{\mathbf{N}} \int_0^t\langle\; \mathcal{R}_{\mathbf{N},\varepsilon}(r)\, ,\,  \big( \nabla_{\varepsilon}^{1,+}\varphi \big)^2+\big( \nabla_{\varepsilon}^{1,-}\varphi \big)^2+\big( \nabla_{\varepsilon}^{2,+}\varphi \big)^2+\big( \nabla_{\varepsilon}^{2,-}\varphi \big)^2 \; \rangle dr.  \nonumber
	\end{eqnarray}
\end{lem}

\begin{pr}  We give the proof for (\ref{ms}), those of (\ref{mi}) and (\ref{mr}) are similar. For all $ r\le t$, we have 
	$ \dis \delta \langle \mathcal{M}_{\mathbf{N},\varepsilon}^S(r),\varphi \rangle  
	=\varepsilon^2  \sum_{i=1}^{\varepsilon^{-2}} \delta \mathcal{S}_{\mathbf{N},\varepsilon}(r,x_i)\varphi(x_i).$
	Since  for $ y_i \neq  x_i  \pm \varepsilon e_j	$, $ \dis \bigg(\delta  \mathcal{S}_{\mathbf{N},\varepsilon}(r,x_i) \bigg)\bigg(\delta  \mathcal{S}_{\mathbf{N},\varepsilon}(r, y_i) \bigg) = 0 $, so
	\begin{eqnarray}\label{saut}
	\Big(\delta \langle \mathcal{M}_{\mathbf{N},\varepsilon}^S(r),\varphi \rangle  \Big)^2
	&=& \varepsilon^4\sum_{i=1}^{\varepsilon^{-2}} \bigg(\delta \mathcal{S}_{\mathbf{N},\varepsilon}(r,x_i)\bigg)^2\varphi^2(x_i) \\
	&+& 2 \varepsilon^4\sum_{i=1}^{\varepsilon^{-2}} \bigg(\delta \mathcal{S}_{\mathbf{N},\varepsilon}(r,x_i) \bigg)\bigg(\delta \mathcal{S}_{\mathbf{N},\varepsilon}(r, x_i+\varepsilon e_1 ) \bigg)\varphi(x_i)\varphi(x_i+\varepsilon e_1) \nonumber \\
	&+& 2 \varepsilon^4\sum_{i=1}^{\varepsilon^{-2}} \bigg(\delta \mathcal{S}_{\mathbf{N},\varepsilon}(r,x_i) \bigg)\big(\delta \mathcal{S}_{\mathbf{N},\varepsilon}(r, x_i-\varepsilon e_1 ) \bigg)\varphi(x_i)\varphi(x_i-\varepsilon e_1)\nonumber \\
	&+& 2 \varepsilon^4\sum_{i=1}^{\varepsilon^{-2}} \bigg(\delta \mathcal{S}_{\mathbf{N},\varepsilon}(r,x_i) \bigg)\bigg(\delta \mathcal{S}_{\mathbf{N},\varepsilon}(r, x_i+\varepsilon e_2 ) \bigg)\varphi(x_i)\varphi(x_i+\varepsilon e_2) \nonumber \\
	&+& 2 \varepsilon^4\sum_{i=1}^{\varepsilon^{-2}} \bigg(\delta \mathcal{S}_{\mathbf{N},\varepsilon}(r,x_i) \bigg)\bigg(\delta \mathcal{S}_{\mathbf{N},\varepsilon}(r, x_i-\varepsilon e_2 ) \bigg)\varphi(x_i)\varphi(x_i-\varepsilon e_2)\nonumber.  
	\end{eqnarray}
	Using successively (\ref{s}) and (\ref{sd}) from the previous lemma, we obtain 
	
	\vspace{-0.8cm}
	
	\begin{eqnarray}\label{cs}
	\sum_{0\le r\le t }\bigg(\delta  \mathcal{S}_{\mathbf{N},\varepsilon}(r, x_i) \bigg)^2 \varphi^2(x_i)
	&\!\!=\!\!&\dfrac{1}{\mathbf{N}}\int_0^t\dfrac{\beta(x_i)\mathcal{S}_{\mathbf{N},\varepsilon}(r,x_i)\mathcal{I}_{\mathbf{N},\varepsilon}(r,x_i)}{\mathcal{S}_{\mathbf{N},\varepsilon}(r,x_i)+\mathcal{I}_{\mathbf{\mathbf{N}},\varepsilon}(r,x_i)+\mathcal{R}_{\mathbf{\mathbf{N}},\varepsilon}(r,x_i)}\varphi^2(x_i)dr \\
	&&\hspace{-8cm}+ \dfrac{4\mu_S}{\mathbf{N}\varepsilon^2} \int_0^t\mathcal{S}_{\mathbf{N},\varepsilon}(r,x_i)\varphi^2(x_i)dr+\dfrac{\mu_S}{\mathbf{N}\varepsilon^2}\int_0^t \bigg[\sum_{j=1}^{2} \mathcal{S}_{\mathbf{N},\varepsilon}(r, x_i  +\varepsilon e_j) + \sum_{j=1}^{2} \mathcal{S}_{\mathbf{N},\varepsilon}(r, x_i-\varepsilon e_j) \bigg] \varphi^2(x_i)dr + \text{Martingale}\nonumber
	\end{eqnarray}
	and
	\begin{eqnarray}\label{sc}
	\sum_{0\le r\le t }\bigg(\delta  \mathcal{S}_{\mathbf{N},\varepsilon}(r,x_i) \bigg)\bigg(\delta  \mathcal{S}_{\mathbf{N},\varepsilon}(r, x_i \pm \varepsilon e_j)\bigg)\varphi(x_i)\varphi(x_i  \pm \varepsilon e_j )
	&& \\
	&&\hspace{-9cm}=-\dfrac{\mu_S}{\mathbf{N}\varepsilon^2} \int_0^t\Big( \mathcal{S}_{\mathbf{N},\varepsilon}(r,x_i)+\mathcal{S}_{\mathbf{N},\varepsilon}(r, x_i \pm \varepsilon e_j)\Big)\varphi(x_i)\varphi(x_i  \pm \varepsilon e_j )dr + \text{Martingale}\nonumber. 
	\end{eqnarray}
	Combining (\ref{saut}), (\ref{cs}) and (\ref{sc}), we deduce that
	
	\begin{eqnarray}
	\sum_{0\le r\le t}\Big(\delta \langle \; \mathcal{M}_{\mathbf{N},\varepsilon}^S(r),\varphi \; \rangle  \Big)^2
	&=& \dfrac{\varepsilon^2}{\mathbf{N}} \int_0^t\langle \; \dfrac{\beta_{\ep}(.)\mathcal{S}_{\mathbf{N},\varepsilon}(r)\mathcal{I}_{\mathbf{N},\varepsilon}(r)}{\mathcal{S}_{\mathbf{N},\varepsilon}(r)+\mathcal{I}_{\mathbf{\mathbf{N}},\varepsilon}(r)+\mathcal{R}_{\mathbf{\mathbf{N}},\varepsilon}(r)}\, , \, \varphi^2 \; \rangle \, dr \nonumber\\
	&+& \dfrac{\mu_S}{\mathbf{N}} \int_0^t\langle \; 4\mathcal{S}_{\mathbf{N},\varepsilon}(r)+ \sum_{j=1}^{2} \mathcal{S}_{\mathbf{N},\varepsilon}(r, .  +\varepsilon e_j) + \sum_{j=1}^{2} \mathcal{S}_{\mathbf{N},\varepsilon}(r, . -\varepsilon e_j) \, , \,  \varphi^2 \; \rangle \, dr\nonumber \\
	&&\hspace{-1.3cm}-\dfrac{2\mu_S}{\mathbf{N}} \int_0^t\langle \; \mathcal{S}_{\mathbf{N},\varepsilon}(r)\, , \,  \sum_{j=1}^{2} \varphi(.)\varphi( . +  \varepsilon e_j ) + \sum_{j=1}^{2} \varphi(.)\varphi( . - \varepsilon e_j )\; \rangle \, dr + \text{Martingale}, \nonumber 
	\end{eqnarray}
	
	which can also be written as
	
	\vspace{-0.8cm}
	
	\begin{eqnarray}
	\sum_{0\le r\le t}\big(\delta \langle \; \mathcal{M}_{\mathbf{N},\varepsilon}^S(r) \, , \, \varphi \, \rangle \big)^2 &=&\dfrac{\varepsilon^2}{\mathbf{N}} \int_0^t\langle \;\dfrac{\beta_{\ep}(.)\mathcal{S}_{\mathbf{N},\varepsilon}(r)\mathcal{I}_{\mathbf{N},\varepsilon}(r)}{\mathcal{S}_{\mathbf{N},\varepsilon}(r)+\mathcal{I}_{\mathbf{\mathbf{N}},\varepsilon}(r)+\mathcal{R}_{\mathbf{\mathbf{N}},\varepsilon}(r)} \, , \,  \varphi^2  \; \rangle \, dr \nonumber \\
	&&\hspace{-3cm} -\dfrac{\mu_S\,\varepsilon^2}{\mathbf{N}} \int_0^t\langle\; \mathcal{S}_{\mathbf{N},\varepsilon}(r)\, ,\,  \big( \nabla_{\varepsilon}^{1,+}\varphi \big)^2+\big( \nabla_{\varepsilon}^{1,-}\varphi \big)^2+\big( \nabla_{\varepsilon}^{2,+}\varphi \big)^2+\big( \nabla_{\varepsilon}^{2,-}\varphi \big)^2 \; \rangle \, dr  + \text{Martingale}. \nonumber
	\end{eqnarray}
	\begin{flushright}
		$\blacksquare$
	\end{flushright}
\end{pr}
The following Lemma generalizes \ref{mmj} in the case of a non constant $\varphi \in C\big(\mathbb{R}_+ ; H^{\varepsilon}\big)$.

\begin{lem}\label{m2}
	The assertion of \ref{mmj}  is valid if $\varphi \in C\big(\mathbb{R}_+ ; H^{\varepsilon}\big)$.
\end{lem}
\nn\begin{pr}
	The general result follows by approximation. $\varphi$ being continuous with respect to $t$, there exists a sequence   $ (\varphi_j)_{_{1\le j\le n}}$ of step functions which converges to $\varphi$ locally uniformly in $[0,\infty). $ It then suffices to consider the case where  $\varphi $ is a step function which we assume from now on. There exists a sequence $0=t_0< t_1<t_2 < \cdots < t_n = t $ such that $\dis \varphi(t,x_i)= \sum_{j=1}^{n} \varphi_j(x_i)\mathbf{1}_{(t_{j-1}, t_j]}(t)$, where $ \varphi_j\in H^{\varepsilon}$, for all $j = 1,\cdots, n$. Applying \ref{mmj} on each interval $(t_{j-1},t_j]$ and summing for all  $j\in\{1,\cdots,n\}$ yields to the result.
	\fpr
\end{pr}
\nn Now we are in a position to give the \\
\textbf{Proof of \ref{ny}} 
\vspace{0.1cm}

Let us fix $ \bar{t} \in (0,T]$ ,   $ \dis i \in \big\{ \;  1, \cdots , \varepsilon^{-2} \; \big\}$ and we use the notation 
$\dis f= \varepsilon^{-2}\mathbf{1}_{V_i}$. We define \\ $ \dis \overline{m}_{\mathbf{N},\varepsilon}^S(t) := \big\langle \int_0^t \mteps(\bar{t}-~r)d\mathcal{M}_{\mathbf{N},\varepsilon}^S(r )\, , f \; \big\rangle  , $ \;    $ 0 \le t \le \bar{t}$. Note that the process $ \dis \left\{ \; \overline{m}_{\mathbf{N},\varepsilon}^S(t) , \;  t \in [ 0,\bar{t} \;] \; \right\} $ is a mean zero martingale and we have\\ $ \dis \overline{m}_{\mathbf{N},\varepsilon}^S(\overline{t}) =  Y_{\mathbf{N},\varepsilon}^S\big(\overline{t}\big) $. 
We have  $ \dis \sum_{0\le r\le t }\Big( \delta  \; \overline{m}_{\mathbf{N},\varepsilon}^S(r)  \, \Big)^2 =  \sum_{0\le r\le t }\bigg( \big\langle \; \delta\mathcal{M}_{\mathbf{N},\varepsilon}^S(r), \mteps(\bar{t}-r)f \; \big\rangle \bigg )^2. $ From \ref{m2}, we have that
\begin{eqnarray}
\sum_{0\le r\le t }\Big( \delta  \; \overline{m}_{\mathbf{N},\varepsilon}^S(r)  \, \Big)^2 -\int_0^tg_{\varepsilon}(r)dr
\end{eqnarray}
is a  mean zero martingale, where
\vspace{-01cm}
\begin{eqnarray}
g_{\varepsilon}(r)
\!\!\!\!&=&\!\!\! \dfrac{\varepsilon^2}{\mathbf{N}}\langle \; \dfrac{\beta(.)\mathcal{S}_{\mathbf{N},\varepsilon}(r)\mathcal{I}_{\mathbf{N},\varepsilon}(r)}{\mathcal{S}_{\mathbf{N},\varepsilon}(r)+\mathcal{I}_{\mathbf{\mathbf{N}},\varepsilon}(r)+\mathcal{R}_{\mathbf{\mathbf{N}},\varepsilon}(r)}\, , \,  \big(\mteps(\bar{t}-r)f\big)^2  \; \rangle  \nonumber \\
&&\hspace{-1.3cm} +\dfrac{\mu_S\, \varepsilon^2}{\mathbf{N}} \langle\; \mathcal{S}_{\mathbf{N},\varepsilon}(r)\, ,\,  \big( \nabla_{\varepsilon}^{1,+}\mteps(\bar{t}-r)f \big)^2+\big( \nabla_{\varepsilon}^{1,-}\mteps(\bar{t}-r)f \big)^2+\big( \nabla_{\varepsilon}^{2,+}\mteps(\bar{t}-r)f \big)^2+\big( \nabla_{\varepsilon}^{2,-}\mteps(\bar{t}-r)f \big)^2 \; \rangle.  \nonumber
\end{eqnarray}
We have
\begin{eqnarray}
g_{\varepsilon}(r)
\!\!\!\!&\le& \!\!\!\! \dfrac{\bar{\beta}\varepsilon^2}{\mathbf{N}}\langle \; 1\, , \,  \big(\mteps(\bar{t}-r)f\big)^2  \; \rangle\nonumber\\
&&\hspace{-1.2cm} +\dfrac{C\mu_S\, \varepsilon^2}{\mathbf{N}} \langle\; 1 \, ,\,  \big( \nabla_{\varepsilon}^{1,+}\mteps(\bar{t}-r)f \big)^2+\big( \nabla_{\varepsilon}^{1,-}\mteps(\bar{t}-r)f \big)^2+\big( \nabla_{\varepsilon}^{2,+}\mteps(\bar{t}-r)f \big)^2+\big( \nabla_{\varepsilon}^{2,-}\mteps(\bar{t}-r)f \big)^2 \; \rangle.  \nonumber
\end{eqnarray}
For $ \theta \in [0,1] $, we define $ m_{\mathbf{N},\varepsilon}^S(t)= \theta \, \mathbf{N}  \, \overline{m}_{\mathbf{N},\varepsilon}^S(t)$. 
$ m_{\mathbf{N},\varepsilon}^S $ is a  mean zero martingale. Furthermore  
\begin{eqnarray}
\vert \delta m_{\mathbf{N},\varepsilon}^S \vert  &\le& \theta \, \mathbf{N}\Big\Vert \mteps(\bar{t}-t)\delta\mathcal{M}_{\mathbf{N},\ep}(t)\Big\Vert_{\infty}\int_Df(x)dx \nonumber\\
&\le& 1 . \nonumber
\end{eqnarray} 
It  follows from \ref{dm} and \ref{em}  that
\begin{eqnarray}
\E\Big( \exp(m_{\mathbf{N},\varepsilon}^S(\bar{t}))\Big) 
& \le & \exp \Big[\dfrac{e}{2}\theta^2 C(\bar{\beta},\mu_S)\mathbf{N} \varepsilon^2(\bar{t}+C\varepsilon^{-2}) \Big] .\nonumber
\end{eqnarray}
It then follows that for any site $x_i \in D_{\varepsilon}$, $\eta > 0 $ 
\begin{eqnarray}
\P\Big( Y_{\mathbf{N},\varepsilon}^S\big(\bar{t}, x_i\big)> \eta \Big)
& =& \P\Big(\theta \, \mathbf{N} \, Y_{\mathbf{N},\varepsilon}^S\big(\bar{t},x_i\big)>\theta \, \mathbf{N} \, \eta \Big) \nonumber \\
& \le & \E\Big[ \exp\Big(\theta \, \mathbf{N} \, Y_{\mathbf{N},\varepsilon}^S\big(\bar{t},x_i\big)\Big) \Big]\exp\big({- \theta \, \mathbf{N} \, \eta }\big) \nonumber \\
& \le & \exp \bigg[\theta \mathbf{N}\Big(C(T) \theta - \eta\Big)\bigg] \quad \text{with} \; C(T)= \dfrac{e}{2}C(\bar{\beta}, \mu_S)(T+C) . \nonumber
\end{eqnarray}
The optimal $\theta $ is $ \theta = \dfrac{\eta }{2C(T)}$, hence
$\dis \P\Big( Y_{\mathbf{N},\varepsilon}^S\big(\bar{t},x_i\big)> \eta \Big) \le \exp\big(-a\, \eta^2\mathbf{N}\big),\; \text{with} \;  a =  \dfrac{1}{4C(T)}. $
We can make a similar computation for $ \P\Big( -Y_{\mathbf{N},\varepsilon}^S\big(\bar{t},x_i\big)> \eta \Big) $ to show that 
$\dis \P\Big( -Y_{\mathbf{N},\varepsilon}^S\big(\overline{t},x_i\big) > \eta \Big) \le \exp(-a \,\eta^2\mathbf{N}). $

Hence for all $ t\in [0,T] $ and $ i \in \left\{ 1 , \cdots , \varepsilon^{-2} \right\} $, we have $$ \P\Big( \big\vert Y_{\mathbf{N},\varepsilon}^S\big(t,x_i\big)\big\vert > \eta \Big) \le 2\exp(-a \eta^2\mathbf{N}).$$
Since $ \dis \Big\Vert Y_{\mathbf{N},\varepsilon}^S(t)\Big\Vert_{\infty} = \un{i}{\text{sup}}\Big\vert Y_{\mathbf{N},\varepsilon}^S\big(t,x_i\big)\Big\vert $,
\begin{eqnarray}\label{yt}
\P\Bigg(\Big\Vert Y_{\mathbf{N},\varepsilon}^S(t)\Big\Vert_{\infty}> \eta \Bigg)
& \le & \sum_{i=1}^{\varepsilon^{-2}}\P\Bigg( \Big\vert Y_{\mathbf{N},\varepsilon}^S\big(t,x_i\big)\Big\vert > \eta \Bigg) \nonumber \\
& \le & 2\varepsilon^{-2} \exp(-a \eta^2\mathbf{N} ).\label{mb}
\end{eqnarray}
We now show that an inequality similar to  (\ref{mb}) holds with  $ \Big\Vert Y_{\mathbf{N},\varepsilon}^S(t)\Big\Vert_{\infty}$ replaced by $ \un{t\in [0,T]}{\sup}\Big\Vert Y_{\mathbf{N},\varepsilon}^S(t)\Big\Vert_{\infty} $ . \\
To this end , we divide $[0,T]$ into $ \varepsilon^{-2} $  intervals $ [nT\varepsilon^2 , (n+1)T\varepsilon^2]$, $0\le n\le \varepsilon^{-2}-1$. 

For  $ t \in [nT\varepsilon^2 , (n+1)T\varepsilon^2]  $, we have \\

\vspace{-0.4cm}

$\dis Y_{\mathbf{N},\varepsilon}^S(t ) =  Y_{\mathbf{N},\varepsilon}^S(nT\varepsilon^2) + \int_{nT\varepsilon^2}^t \Delta_{\varepsilon} Y_{\mathbf{N},\varepsilon}^S(r)dr + \tilde{m}_{\mathbf{N},\varepsilon}^S (t)$, 
where $ \dis \tilde{m}_{\mathbf{N},\varepsilon}^S (t) = \mathcal{M}_{\mathbf{N},\varepsilon}^S(t)-  \mathcal{M}_{\mathbf{N},\varepsilon}^S(nT\varepsilon^2). $  \\
We have 
\begin{eqnarray} \label{ysup}
\hspace{-1.5cm}\Big\Vert Y_{\mathbf{N},\varepsilon}^S(t) \Big\Vert_{\infty} 
& \le &  \Big\Vert Y_{\mathbf{N},\varepsilon}^S(nT\varepsilon^2)\Big\Vert_{\infty} + 8 \; \varepsilon^{-2} \int_{nT\varepsilon^2}^t  \Big\Vert Y_{\mathbf{N},\varepsilon}^S(r) \Big\Vert_{\infty} dr + \Big\Vert \tilde{m}_{\mathbf{N},\varepsilon}^S (t) \Big\Vert_{\infty},
\end{eqnarray}
so  Gronwall's inequality implies that

\vspace{-1.1cm}

\begin{eqnarray}\label{sg}
\un{t \in [nT\varepsilon^2 , (n+1)T\varepsilon^2]}{\text{sup}}\Big\Vert Y_{\mathbf{N},\varepsilon}^S(t) \Big\Vert_{\infty} 
& \leq & \Bigg( \Big\Vert Y_{\mathbf{N},\varepsilon}^S(nT\varepsilon^2)\Big\Vert_{\infty} + \un{t \in [nT\varepsilon^2 , (n+1)T\varepsilon^2]}{\sup}\Big\Vert \tilde{m}_{\mathbf{N},\varepsilon}^S (t) \Big\Vert_{\infty} \Bigg)\exp(8T) .
\end{eqnarray}
We now fix $ i \in \left\{  1 , \cdots , \varepsilon^{-2}  \right\}  , \theta \in [0,1] $ and  set $ m_{\mathbf{N},\varepsilon}^S(t)= \theta \, \mathbf{N}\,\tilde{m}_{\mathbf{N},\varepsilon}^S \big(t\big). $ 
It follows from \ref{mart1} that
\begin{align} 
\sum_{ nT\varepsilon^2 \le r\ \le t } \Big(\delta m_{\mathbf{N},\varepsilon}^S(r)\Big)^2 - \dfrac{\mu_S}{\varepsilon^2}\theta^2\mathbf{N}\int_{nT\varepsilon^2}^t\Big(\sum_{y\sim x_i}\mathcal{S}_{\mathbf{N},\varepsilon}(r,y) + 4 \mathcal{S}_{\mathbf{N},\varepsilon}(r,x_i)\Big)dr 
& \\
& \hspace{-10cm}-\theta^2 \mathbf{N} \int_{nT\varepsilon^2}^t \dfrac{\beta(x_i)\mathcal{S}_{\mathbf{N},\varepsilon}(r,x_i)\mathcal{I}_{\mathbf{N},\varepsilon}(r,x_i)}{\mathcal{S}_{\mathbf{N},\varepsilon}(r,x_i)+\mathcal{I}_{\mathbf{\mathbf{N}},\varepsilon}(r,x_i)+\mathcal{R}_{\mathbf{\mathbf{N}},\varepsilon}(r,x_i)}dr \nonumber
\end{align}
is a  mean zero martingale and $ \Big\vert \delta  m_{\mathbf{N},\varepsilon}^S(t) \Big\vert  \le 1 $ .  
Furthermore, for $ nT\varepsilon^2 < t \le (n+1)T\varepsilon^2$
\begin{eqnarray}
\hspace{-1cm}\dfrac{\mu_S}{\varepsilon^2}\theta^2\mathbf{N} \int_{nT\varepsilon^2}^t\Big(\sum_{y\sim x_i}\mathcal{S}_{\mathbf{N},\varepsilon}(r,y) + 4 \mathcal{S}_{\mathbf{N},\varepsilon}(r,x_i)\Big) dr +  \mathbf{N}\theta^2 \int_{nT\varepsilon^2}^t \dfrac{\beta(x_i)\mathcal{S}_{\mathbf{N},\varepsilon}(r,x_i)\mathcal{I}_{\mathbf{N},\varepsilon}(r,x_i)}{\mathcal{S}_{\mathbf{N},\varepsilon}(r,x_i)+\mathcal{I}_{\mathbf{\mathbf{N}},\varepsilon}(r)+\mathcal{R}_{\mathbf{\mathbf{N}},\varepsilon}(r,x_i)} dr\nonumber\\
&& \hspace{-3cm}\le C(\bar{\beta},\mu_S)T \mathbf{N} \theta^2. \nonumber 
\end{eqnarray}

Hence  by  \ref{em}, it follows that
$\dis  \E\Big[ \exp\Big(m_{\mathbf{N},\varepsilon}^S\big((n+1)T\varepsilon^2 \big) \Big)\Big] \le \exp\Big[ C(\bar{\beta}, \mu_S) \mathbf{N} \theta^2 T\Big].$

It then follows from Doob's inequality that
\begin{eqnarray}
\P\Bigg( \un{t \in [nT\varepsilon^2 , (n+1)T\varepsilon^2]}{\text{sup}} \tilde{m}_{\mathbf{N},\varepsilon}^S \big(t,x_i\big)  \ge \eta \Bigg)
&\le & \E\Big[ \exp\Big(m_{\mathbf{N},\varepsilon}^S\big((n+1)T\varepsilon^2 \big) \Big)\Big] \exp(-\theta \mathbf{N}\eta) \nonumber \\
&\le & \exp\Big[ \theta \mathbf{N}\Big(C(T)\,\theta - \eta \Big) \Big].\nonumber 
\end{eqnarray}
Choosing $ \theta = \dfrac{\eta}{2C(T)}$, we deduce that 

$$ \P\Bigg( \un{t \in [nT\varepsilon^2 , (n+1)T\varepsilon^2]}{\text{sup}} \tilde{m}_{\mathbf{N},\varepsilon}^S \big(t,x_i\big)  \ge \eta \Bigg) \le 
\exp(-a \eta^2\mathbf{N}),\; \;  \text{where}\; \;  a  = \dfrac{1}{4C(T)}. $$ 
The same hold for $ - \tilde{m}_{\mathbf{N},\varepsilon}^S \big(t,x_i\big) $. Consequently
\begin{eqnarray}\label{msup}
\P\Bigg( \un{t \in [nT\varepsilon^2 , (n+1)T\varepsilon^2]}{\sup}\Big\Vert \tilde{m}_{\mathbf{N},\varepsilon}^S \big(t,x_i\big) \Big\Vert_{\infty} \ge \eta \Bigg) 
& \le & 2 \; \varepsilon^{-2} \exp(-a \eta^2\mathbf{N}).
\end{eqnarray} 
Combining the inequalities (\ref{yt}), (\ref{sg}) and (\ref{msup}), we obtain
\begin{eqnarray}
\P\Bigg( e^{-8T} \un{t \in [nT\varepsilon^2 , (n+1)T\varepsilon^2]}{\text{sup}}\Big\Vert Y_{\mathbf{N},\varepsilon}^S(t) \Big\Vert_{\infty} \ge \eta \Bigg) 
& \le & 4 \;   \varepsilon^{-2} \exp(-a \dfrac{\eta^2}{4} \mathbf{N}) , 
\end{eqnarray}
from which we deduce that
\begin{eqnarray}\label{martsup}
\P\Bigg( e^{-8T}\un{t \in [0,T]}{\text{sup}}\Big\Vert Y_{\mathbf{N},\varepsilon}^S (t) \Big\Vert_{\infty} \ge \eta \Bigg) 
& \le & \sum_{n=0}^{\varepsilon^{-2}-1} \P\Bigg( e^{-8T}\un{t \in [nT\varepsilon^2 , (n+1)T\varepsilon^2]}{\text{sup}}\Big\Vert Y_{\mathbf{N},\varepsilon}^S (t) \Big\Vert_{\infty} \ge \eta \Bigg) \nonumber \\
& \le & 4 \; \varepsilon^{-4} \exp(-a \dfrac{\eta^2}{4} \mathbf{N}). 
\end{eqnarray}
Since $ \dfrac{\mathbf{N}}{\log (1/\varepsilon)} \longrightarrow +\infty $ implies that  $ \varepsilon^{-4} \exp(-a \eta^2 \mathbf{N}) \longrightarrow 0 $, we have proved that  $ \dis \un{t \in [0,T]}{\sup}\Big\Vert Y_{\mathbf{N},\varepsilon}^S (t) \Big\Vert_{\infty} \longrightarrow 0 $ in probability. The same arguments show that $ \dis \un{t \in [0,T]}{\sup}\Big\Vert Y_{\mathbf{N},\varepsilon}^I (t) \Big\Vert_{\infty}+ \un{t \in [0,T]}{\sup}\Big\Vert Y_{\mathbf{N},\varepsilon}^R(t) \Big\Vert_{\infty} \longrightarrow 0 $ in probability  as $\mathbf{N}\rightarrow \infty$ and $ \varepsilon \to 0 $, under our standing assumption. Finally, we have shown that $  \dis \un{t \in [0,T]}{\sup}\Big\Vert Y_{\mathbf{N},\varepsilon} (t) \Big\Vert_{\infty} \longrightarrow 0 $ in probability, which completes the proof of the Proposition. 

\vspace{-0.3cm}

\begin{flushright}
	$\blacksquare$
\end{flushright}

\begin{Large}\textbf{Remark 1}\end{Large}


The law of large numbers in sup-norm remains true in dimensions $d=1,3$. To  see that, it suffices to remark that  $ \dis \Delta_{\varepsilon}= \sum_{j=1}^{d}\nabla_{\varepsilon}^{j,-}\nabla_{\varepsilon}^{j,+} $ has always $\varepsilon^{-d}$ bounded eigenvectors. In this case the \ref{dm} become 

$$ \Big\langle \;  \sum_{j=1}^{d}\Big(\nabla_{\varepsilon}^{j,+}\mtepj(t)f\Big)^2 + \sum_{j=1}^{d}\Big(\nabla_{\varepsilon}^{j,-}\mtepj(t)f\Big)^2 + \Big(\mtepj(t)f\Big)^2 , 1 \;  \Big\rangle  \le h_{\varepsilon}(t) $$
where \; $\dis \int_0^t h_{\varepsilon}(r)dr\le C\,\varepsilon^{-d} + t. $
Hence (\ref{martsup}) becomes \; 
$\dis
\P\Bigg( e^{-8T}\un{t \in [0,T]}{\text{sup}}\Big\Vert Y_{\mathbf{N},\varepsilon}^S (t) \Big\Vert_{\infty} \ge \eta \Bigg) 
\le 4 \; \varepsilon^{-d-2} \exp(-a \dfrac{\eta^2}{4} \mathbf{N}). $ \\

Moreover, the result holds for periodic boundary conditions. Indeed, in this case, the eigenvectors of the Laplace operator are the product of the one-dimensional eigenvectors

$$\varphi_{n}(x)=\left\{
\begin{array}{cl}
& 1 , \; \mbox{for} \; \;  n= 0 ,\\
& \sqrt{2}\cos(n\pi  x), \; \text{for} \; \; n>0 \; \text{and even},
\end{array}
\right.
$$
$\hspace{4.7cm}\dis \psi_{n}(x)=\sqrt{2}sin(n\pi x), \; \text{for} \; \; n>0 \; \text{and even}. $ \\

\begin{Large}\textbf{Remark 2}\end{Large}

We conclude that, by two laws of large numbers, the consistency of the various  models has been established. \\
In a furture  work, we will study the fluctuations of the  stochastic model around its  deterministic law of large numbers limit.\\

\textbf{Acknowledgments}.  The authors are deeply indebted to the referee for a careful reading and several suggestions that greatly improved the paper.\\

\end{document}